\documentclass[10pt]{article}
\usepackage{amssymb}
 \newtheorem{prop}{Proposition}
 \newtheorem{Definition}{Definition}
 \newtheorem{theo}{Theorem}
 \newtheorem{lemma}{Lemma}
 
 \newcommand{\bs}{\left\{}
 \newcommand{\es}{\right.}
 \newcommand{\ba}{\begin{array}}
 \newcommand{\ea}{\end{array}}
 \newcommand{\be}{\begin{equation}}
 \newcommand{\ee}{\end{equation}}

 \def\RR{{\rm I\hspace{-0.50ex}R} }

 \def\CC{\rm \hbox{C\kern-.57em\raise.47ex
 \hbox{$\scriptscriptstyle |$}\kern+0.5 em}}
 
 \title{Study of the linear ablation growth rate for the quasi-isobaric model of Euler equations with thermal conductivity}
 \author{Olivier Lafitte\thanks{CEA/DM2S, Centre
 d'Etudes de Saclay, 91191 Gif sur Yvette Cedex}\thanks{Universit{\'e} de Paris XIII, LAGA, 93 430 Villetaneuse}}
 \date{November 30, 2006}

\begin{document}
 \maketitle

 \begin{abstract}
 In this paper, we study a linear system related to the 2d system of Euler equations with thermal conduction in the quasi-isobaric approximation of Kull-Anisimov \cite{Kull}. This model is used for the study of the ablation front instability, which appears in the problem of inertial confinement fusion. The heat flux ${\vec Q}$ is given by the Fourier law $T^{-\nu}{\vec Q}$ proportional to $\nabla T$, where $\nu>1$ is the thermal conduction index, and the external force is a gravity field ${\vec g}=-g{\vec e}_x$. This physical system contains a mixing region, in which the density of the gaz varies quickly, and one denotes by $L_0$ an associated characteristic length. The fluid velocity in the denser region is denoted by $V_a$. \\

The system of equations is linearized around a stationary solution, and each perturbed quantity ${\tilde u}$ is written using the normal modes method
$${\tilde u}(x, z, t)=\Re ({\bar u}(x, k, \gamma)e^{ikz+\gamma\sqrt{gk} t})$$
in order to take into account an increasing solution in time.\\
 The resulting linear system is a non self-adjoint fifth order system. Its coefficients depend on $x$ and on physical parameters $\alpha, \beta$, $\alpha$ and $\beta$ being two dimensionless physical constants, given by $\alpha\beta = kL_0$ and $\frac{\alpha}{\beta}=\frac{gL_0}{V_a^2}$ (introduced in \cite{CCLARA}). We study the existence of bounded solutions of this system in the limit $\alpha\rightarrow 0$, under the condition $\beta \in [\beta_0, \frac{1}{\beta_0}]$, and the assumption $\Re\gamma\in [0, \frac{1}{\beta_0}], |\gamma|\leq \frac{1}{\beta_0}$ (regime that we studied for a simpler model in \cite{CCLARA}) calculating the Evans function $Ev(\alpha, \beta, \gamma)$ associated with this linear system.\\
Using rigorous constructions of decreasing at $\pm\infty$ solutions of systems of ODE, we prove that, for $\beta \in [\beta_0, \frac{1}{\beta_0}]$, $\Re \gamma\in [0, \frac{1}{\beta_0}]$, $|\gamma|\leq\frac{1}{\beta_0}$, there exists $\alpha_1>0$ such that there is no bounded solution of the linearized system for $0<\alpha\leq \alpha_1$. \\Necessarily, for any $M>0$ and $\beta_0>0$ there exists $\alpha_1>0$ such that, for $0<\alpha\leq \alpha_1$ and $\beta\in [\beta_0, \frac{1}{\beta_0}]$, an admissible value $\gamma(\alpha, \beta)$ such that there exists a bounded solution of the linearized system satisfying $|\gamma|\leq M$ is such that $\Re\gamma\notin [0, M]$.\\
\end{abstract}
\tableofcontents
\setcounter{section}{-1}
\section{Introduction}
This paper is devoted to the precise calculus of the Evans function $Ev(\alpha, \beta, \gamma)$ of the normal mode formulation of the linearized system of equations associated with the quasi-isobaric (low Mach number) model. The calculus of this Evans function is not classical, because the matrix of the differential system has singular coefficients,and because these coefficients do not behave exponentially in the spatial variable. However, usual techniques of ordinary differential equations and introduction of a Fuchsian problem allow us to calculate this Evans function under certain assumptions on the parameters $\alpha, \beta, \gamma$ introduced in the Abstract.\\
In this Introduction, we first describe the physical model (\ref{Eulerkull}), define what is called a linear growth rate of the linearized system associated with a stationary solution of this physical model, then finally describe the contents of each step of the proof of the main Theorem (Theorem \ref{lacanau}).

\subsection{Physical model}

We consider a compressible fluid characterized by its density $\rho$, its velocity $(u, v)$ and its temperature $T$ in a gravity field ${\vec g}= -|g|{\vec e}_x$. We assume that this fluid has the following properties:

a) when $x$ goes to $+\infty$, for all $z$ we have $\rho\rightarrow \rho_a$, $(u, v)\rightarrow (-V_a, 0)$ and  $T\rightarrow T_a$.

b) the functions $(\rho, u, v, p, T)$, where $p$ is linked to the pressure in the fluid, satisfy the system of the Euler equations in two dimensions $(x, z)$ with thermal conduction in the quasi-isobaric approximation for a perfect gaz:
\be\label{lowmach}\mbox{div}(C_p\rho T {\vec u}+{\vec Q})=0,\ee
where $C_p$ is the calorific capacity of the fluid, the heat conduuction flux ${\vec Q}$ being given by the Fourier law
\be\label{flux}{\vec Q}=-k(T)\nabla T,\ee
the thermal conduction law is \be\label{conductionthermique}k(T)=K_0T^{\nu},\ee where $\nu$ is the thermal conduction indice.\\
We introduce a characteristic length $L_0$ associated with the thermal properties of the fluid
\be \label{longueur}L_0={K_0T_a^{\nu}\over C_p\rho_aV_a}.\ee
Physical values of $L_0$ for the case of the ICF are of order $10^{-5}$ meters.

Under a quasi-isobaric assumption, the system modelizing the ablation model was given by H.J. Kull \cite{Kull} and appears for example in P. L. Lions \cite{lions}. It states
\be\label{Eulerkull}\bs\ba{l}\partial_t\rho + \mbox{div}(\rho{\vec u})=0\cr
\partial_t(\rho{\vec u})+\mbox{div}(\rho {\vec u}\otimes{\vec u}+p)=\rho{\vec g}\cr
\rho T=\rho_aT_a\cr
{\vec Q}=-k(T)\nabla T\cr
\mbox{div}(C_p\rho T{\vec u}+{\vec Q})=0.\ea\es\ee
This model can be derived either from the low Mach approximation (Majda \cite{Majda}, Dellacherie \cite{Della}) or the quasi-isobaric approximation (Kull \cite{Kull}, Kull-Anisimov \cite{kullA}, Masse \cite{Masse}). See a short analysis in Section \ref{sec1}.
A stationary laminar solution of the system (\ref{Eulerkull}) is $(\rho_0(x), u_0(x), 0, p_0(x), T_0(x))$, where we introduce a function $\xi$ such that
$$\rho_0(x)=\rho_a\xi({x\over L_0}), u_0(x)=-{\rho_aV_a\over \rho_0(x)}=-{V_a\over \xi({x\over L_0})}, T_0(x)={\rho_aT_a\over \rho_0(x)}={T_a\over \xi({x\over L_0})}$$
and $p_0(x)$ satisfies
$$p_0(x)+{\rho_a^2V_a^2\over \rho_0(x)}+g\int_{x_0}^x\rho_0(s)ds= p_0(x_0)+{\rho_a^2V_a^2\over \rho_0(x)}.$$
The function $\xi$ is the solution of the differential equation
\be\label{xi}{d\xi\over dy}=\xi^{\nu+1}(1-\xi)\ee
such that $\xi(0)=\frac{\nu+1}{\nu+2}$. Note that, in this case
$$p_0(L_0y)=p_0(0)+{\nu+2\over \nu+1}\rho_aV_a^2-\rho_aV_a^2[{1\over \xi(y)}+\frac{gL_0}{V_a^2}\int_0^y\xi(t)dt].$$
Introduce
\be Z(\rho)={\rho_a^{\nu+1}\over (\nu+1)\rho^{\nu+1}}.\label{fonctionZ}\ee
The system of unknowns that we consider is
$${\tilde U}=\left(\ba{c}\rho u\cr \rho u^2+p\cr \rho u v\cr Z(\rho)\cr u-L_0V_a\partial_x(Z(\rho))\ea\right).$$
From ${\tilde U}$, we recover $\rho$ from $Z(\rho)$, $u={\rho u \over \rho}$, $v={\rho u v\over \rho u}$ and $p=p+\rho u^2-{(\rho u)^2\over \rho}$.\\
For this choice of unknowns, we introduce $F_1({\tilde U})={\tilde U}$. There exists three explicit functions $F_0$, $F_2$ and $F_3$ such that the system (\ref{Eulerkull}) is equivalent to the system on ${\tilde U}$
\be\label{ssy}\bs\ba{l}\partial_tF_0({\tilde U})+\partial_xF_1({\tilde U})+\partial_z(F_2({\tilde U}, \partial_z{\tilde U}))=F_3({\tilde U})\cr
T={\rho_aT_a\over Z^{-1}({\tilde U}_4)}\ea\es\ee
A stationnary laminar solution $U_0$ of this system depending only on $x$ is $U_0(x)$ such that
$${d\over dx}(F_1(U_0(x)))=F_3(U_0(x)).$$
The identity $F_1({\tilde U})= {\tilde U}$ is the natural choice when one studies a basic solution depending only on the variable $x$.
\subsection{Definition of a linear growth rate}
We linearize (\ref{ssy}) around $U_0(x)$. Denote by $U$ the unknowns ${\tilde U}-U_0(x)$. The linearized system writes ($\nabla_1$ and $\nabla_2$ denotes the gradient of $F_2$ with respect to the first and second set of variables ${\tilde U}$ and $\partial_z{\tilde U}$):
\be\partial_xU+\nabla F_0(U_0(x))\partial_t U+\nabla_1F_2(U_0(x), 0)\partial_zU+\nabla_2F_2(U_0(x), 0)\partial^2_{z^2}U=\nabla F_3(U_0(x))U\label{sl}\ee
which can be rewritten $M(x, \partial_x, \partial_z, \partial^2_{z^2}, \partial_t)U=0$. Note that its coefficients depend on $x$ through the stationnary solution.

We are now ready to introduce the definition of a linear growth rate for a non linear system around a laminar solution:

\begin{Definition}
\label{defigrowth}
Let $M(x, \partial_x, \partial_z, \partial^2_{z^2},\partial_t)U=0$ be the linearized system

We call a linear growth rate of this system for the wave number $k$ a value of $\sigma$ (depending on $k$) such that $\Re\sigma\geq 0$ and there exists a non-trivial solution $U(x, k, \sigma)$ of the system

\be\label{LS}M(x, {d\over dx}, ik, -k^2,\sigma)U(x, k, \sigma)=0\ee
such that $U$ is bounded and going to 0 when $x$ goes to $\pm\infty$.
The function
$$U(x, k, \sigma )e^{ikz}e^{\sigma t}$$
is called a normal mode solution of the system.

\end{Definition}
The normal mode system associated with (\ref{sl}) is

\be\label{star}{dV\over dx}+\nabla F_0(U_0(x))\sigma V+ik\nabla_1F_2(U_0(x), 0)V-k^2\nabla_2F_2(U_0, 0)V=\nabla F_3(U_0(x))V.\ee
The scope of this paper is to find bounded non trivial solutions of (\ref{star}), and associated values of $\sigma$ if any. If such a solution exist, it will lead to a normal mode solution of the linearized system. Note that, in the set-up we described, different physical parameters appear, namely $k, L_0, V_a, g$. As the classical growth rate of Rayleigh is equal to $(\frac{\rho_2-\rho_1}{\rho_2+\rho_1}gk)^{\frac12}$ for the discontinuity model \cite{Rayleigh}, and as we proved (\cite{CCLARA}, \cite{HelLaf}) that this value was the limit of the growth rate when $kL_0$ goes to zero, we are led to introducing the following quantities
\be \varepsilon = kL_0, \quad Fr={V_a^2\over gL_0}\ee
and \be\alpha=\sqrt{\varepsilon\over Fr}, \beta=\sqrt{\varepsilon Fr}, \gamma=\frac{\sigma}{ \sqrt{gk}}.\ee
The aim of this paper is to study the existence of a growth rate $\gamma$, $\Re\gamma\geq 0$, in the limit $L_0\rightarrow 0$  when the Froude number $Fr$ is of order ${1\over L_0}$, which means that
$$\alpha\rightarrow 0, \beta>0.$$
\paragraph{Remark} Other regimes rely on different assumptions on $\alpha$ and $\beta$: we refer to \cite{HL2} for the results that can be obtained for this model in the high frequency regime $k\rightarrow +\infty$. In this other regime the scaling writes\\
\centerline{$\varepsilon$ large, $\varepsilon^3Fr\leq C'$, where $C'$ is a constant.}
 In the first section, we study the physical origin of the model and derive the properties of the stationary solution, where the associated density profile satisfies:
$$\bs\ba{l}\rho_0(x)\rightarrow 0\mbox{ when }x\rightarrow -\infty\cr
\rho_0(x)-\rho_a\simeq Ce^{-{x\over L_0}}, x\rightarrow +\infty\cr
\rho_0(x)|x|^{1\over \nu}\rightarrow \rho_a({L_0\over \nu})^{1\over \nu}, x\rightarrow -\infty.\ea\es$$
We then derive the linearized system, which is a fifth order differential system whose coefficients depend on $\rho_0(x)$ and are singular when $\rho_0(x)$ go to zero. Note that this is not a classical case for the study of such systems and that this leads to rather tricky methods.\\
In the second section, we recall the general set-up for the calculation of the Evans function of the linearized system, and we give the induced differential systems in $\Lambda^n(K^5)$ for $n=2, 3$. Note that the field $K$ is $\RR$ for real values of $\gamma$ and $K=\CC$ for complex values of $\gamma$. This Evans function is (related to) the vectorial product of the normalized solution in $\Lambda^2(K^5)$ which has the greatest decay when $x\rightarrow +\infty$ and of the normalized solution in $\Lambda^3(K^5)$ which has the greatest decay at $x\rightarrow -\infty$.\\
In the third section, we identify the solutions of the system in $\Lambda^2(K^5)$ deduced from (\ref{LS}) for $x\rightarrow +\infty$. In this region, we use the exponential behavior of the profile to obtain the classical analytic expansion of the normalized solution of the system in $\Lambda^2(K^5)$. There exists $\xi_0\in ]0, 1[$ (corresponding to $y_0\in \mathbb R$ through $\xi(y_0)=\xi_0$) such that this analytic expansion is valid for $\xi({x\over L_0})\geq \xi_0$, that is $x\geq L_0y_0$. Note that, however, the expansion of this solution cannot be obtained by the techniques developed in Zumbrun {\it et al} \cite{Zum}, because the Gap lemma assumptions are not fulfilled.\\
A general feature in the calculation of the Evans function is to obtain an overlapping region of definition between the solution well behaved at $+\infty$ and the solution well behaved at $-\infty$. A first step to achieve this overlap is then to prove that there exists $\alpha_0>0$ and $R>0$ such that, for all $\alpha<\zeta<\frac{1}{R}$ the solution obtained for $x\geq L_0y_0$ can be extended in $[X_*(\alpha, \zeta), L_0y_0]$ where
$$\xi(\frac{X_*(\alpha, \zeta)}{L_0})=\alpha^\frac{1}{\nu}\zeta^{-\frac{1}{\nu}}.$$ This is the aim of the fifth section. The behavior of the solution in the region $[\alpha^{1\over \nu}\zeta_0^{-{1\over \nu}}, \xi_0]$ when $\alpha\rightarrow 0$ is different from the classical analytic expansion in $\alpha$ for $y\in [y_0, +\infty[$ and it is the aim of Sections \ref{suniform} and \ref{sprecisees}.\\
 Once this extension is done, an easy calculus is the calculus of a growth rate associated with the following stationnary solution, characterized by its density profile, for a $\zeta_0$ such that $\zeta_0<\frac{1}{R}$:

\be \rho_*(x)=\bs\ba{l}\rho_a\xi(\frac{x}{L_0}), x\geq X_*(\alpha, \zeta_0)\cr
\rho_a\alpha^\frac{1}{\nu}\zeta_0^{-\frac{1}{\nu}}.\ea\es\ee
This calculus is an improvement of the discontinuity model of Piriz, Sanz and Ibanez \cite{Piriz} and it is the aim of Section \ref{petitmodele}.\\
When the profile is not constant in the region $]-\infty, X_*(\alpha, \zeta_0)]$ (that is for the full model), the system leads to a fuchsian problem in the region $x\rightarrow -\infty$, and we use the hypergeometric equation (see \cite{lafitteX}). The solution of the system in $\Lambda^3(K^5)$ deduced from (\ref{LS}) is identified in any region of the form ${x\over L_0}\in ]-\infty,-{t_0\over \alpha\beta}]$ for every $t_0$, which means that $x\in ]-\infty, -{t_0\over k}]$. The results of the analysis of these solutions is summarized in Theorem \ref{latotale}. \\
The study of the roots of the Evans function is the aim of Section \ref{finalO} and we summarize the method here. From the relation $\xi(y)|y|^{1\over \nu}\rightarrow \nu^{-{1\over
    \nu}}$ when $y\rightarrow -\infty$, we deduce that $-\alpha\beta
X_*(\alpha)\rightarrow {\beta\zeta_0\over \nu}>0$ when
$\alpha\rightarrow 0$. Hence for all $<0\alpha\leq \alpha_0$ there
exists $t_0>0$ such that the regions $]-\infty, -{t_0\over
  \alpha\beta}]$ and $[{X_*(\alpha)\over L_0},y_0]$, where
$\xi(y_0)=\xi_0$, overlap.\\
We then express the Evans function $Ev(\alpha, \beta, \gamma)$ of the
system at a point of $[-\frac{t_0}{\alpha\beta},
\frac{X_*(\alpha)}{L_0}]$. The limit when $\alpha\rightarrow 0$ and
$t_0$ small exists and we write its expression in terms of
$r=\frac{\gamma}{\beta}$, $\beta$ and $t_0>0$. As it does not depend
on $t_0$ we study the limit when $t_0\rightarrow +\infty$, hence
proving that the only positive value of $r$ which is admissible is
$r=1$. We deduce a contradiction, proving that there is no growth rate ( of 
positive real part) for the system. This can be stated as

\paragraph{Theorem}
{\bf 
Let $M$ be given. There exists $\alpha_*>0$ such that, for $0<\alpha<\alpha_*$, $\beta\in [\frac{1}{M}, M]$, the Evans function $Ev(\alpha, \beta, \gamma)$ of the system  has no root for $|\gamma|\leq M$, $\Re \gamma\in [0, M]$.}

\section{Derivation of the quasi-isobaric model}
\label{sec1}
\subsection{The physical approximations}
The general equations are the thermal hydrodynamic equations, written in a non conservative form:

\be\label{laphysique}\bs\ba{l}\partial_t\rho+\mbox{div}(\rho{\vec u})=0\cr
\partial_t(\rho{\vec u})+\mbox{div}(\rho{\vec u}\otimes {\vec u}+p)=\rho{\vec g}\cr
\rho(\partial_t+{\vec u}.\nabla )h - (\partial_t+{\vec u}.\nabla )p=-\mbox{div}({\vec Q}+{\vec J}_l).\ea\es\ee
where $C_p$ and $C_v$ are the classical tehrmodynamic calorific capacities at constant pressure and at constant volume, $h$ is the enthalpy $h=C_pT$, the pressure and the density being given by the equation of state $p=(C_p-C_v)\rho T$, ${\vec Q}=-k(T) \nabla T$, ${\vec J}_l=0$ (in our assumption the energy given to the system is 0).
The quasi-isobaric approximation writes

$${C_p-C_v\over C_v}{\delta p\over p}<<{\delta T\over T}.$$
Following L. Masse \cite{Masse}, this relies on the two hypotheses

\be M^2<<1, {M^2\over Fr}<<1.\ee
Hence the quasi-isobaric model relies on a low Mach hypothesis.\\
The equation of the energy rewrites
 $$\mbox{div}{\vec u}+{1\over \rho h}\mbox{div}{\vec Q}= -{(\partial_t+{\vec u}.\nabla)(\rho h -p)\over \rho h}.$$
 This equation is {\bf approximated} by (\ref{lowmach}) as we will see below.
What follows is a {\bf formal} derivation of the quasi-isobaric model under a low mach hypothesis. It is closely related to the method given by Majda \cite{Majda}, Dellacherie \cite{Della}.
\subsection{Adimensionnalization and low Mach expansion}
Use the reference density $\rho_a$, the reference velocity $V_a$, and the reference pressure $p_s$ associated with the sound velocity $c_s$ such that $c_s^2\frac{C_p}{C_v} \rho_a=p_s$. The Mach number is thus $M={V_a\over c_s}$.\\
Write $\rho=\rho_a\rho'$, ${\vec u}=V_a{\vec u}'$, $p=p_sp'$. The system of equations (\ref{laphysique}) rewrites

$$\bs\ba{l}V_a^{-1}\partial_t\rho'+\mbox{div}(\rho'{\vec u}')=0\cr
V_a^{-1}\partial_t(\rho'{\vec u}')+\mbox{div}(\rho'{\vec u}'\otimes {\vec u}'+\gamma{p'\over M^2})={g\over V_a^2}\rho'\cr
\mbox{div}{\vec u}'+V_a^{-1}{1\over C_p\rho T}\mbox{div}{\vec Q}=-{(V_a^{-1}\partial_t+{\vec u}'.\nabla)p\over p}\ea\es$$
If we assume that all the quantities have an asymptotic expansion in $M$, in particular $p'=p'_0(x, z, t)+M^2p(x, z, t, M)$ we have the following relations from the momentum equations

$$\partial_xp'_0=0, \partial_zp'_0=0$$
hence $p'_0$ depends only on $t$. This is the same result as in the analysis of Dellacherie \cite{Della}.
In this model, we assume that the pressure $p'_0$ is constant, because we assume that the ground state for the equations is stationnary.

Replacing the relation $p'(x, z, t, M)=p'_0+M^2p$ in the energy equation we obtain

$$\mbox{div}{\vec u}'+V_a^{-1}{C_p-C_v\over C_p(p'_0+M^2p)}\mbox{div}{\vec Q}=-M^2{(V_a^{-1}\partial_t+{\vec u}'.\nabla)p\over p'_0+M^2p}$$
Finally, using ${\vec Q}=K_0( {p'_0+M^2p\over (C_p-C_v)\rho'})^{\nu}\nabla{p'_0+M^2p\over (C_p-C_v)\rho'}$, we deduce that

$${\vec Q}=K_0({p'_0\over C_p-C_v})^{\nu+1}\nabla Z(\rho)+O(M^2)$$
hence the formal analysis leads to the equation
$$\mbox{div}(V_a^{-1}{\vec u}'+L_0\nabla Z(\rho))=O(M^2)$$
where we used ${p'_0\over C_p-C_v}={\rho_aT_a\over p_s}$ deduced from the relation $p= (C_p-C_v)\rho T$.
The resulting equation can be written $\mbox{div}(C_p\rho T {\vec u}+{\vec Q})=0$, $\rho T=\rho_aT_a$ hence (\ref{lowmach}).\\
Finally, in the momentum equations, rewriting ${p'\over M^2}={p'_0\over M^2}+p$ and using $p'_0$ constant, we obtain the equations
$$V_a^{-1}\partial_t(\rho'{\vec u}')+\mbox{div}(\rho'{\vec u}'\otimes {\vec u}'+p Id)= {g\over L_a^2}\rho'.$$
Note finally that the relation (\ref{lowmach}) and the relation $\rho T= \rho_aT_a$  lead to the equation on $\rho$:
\be\label{equationZ}(\partial_t+{\vec u}.\nabla)Z(\rho) - (\nu+1)L_0V_aZ(\rho)\Delta  Z(\rho)=0.\ee

\subsection{Study of the stationnary solution}
The resulting system of equations that models our phenomenon is thus
\be\label{KANL}\bs\ba{l}\partial_t\rho+\mbox{div}(\rho {\vec u})=0\cr
\partial_t(\rho {\vec u})+\mbox{div}(\rho {\vec u}\otimes {\vec u} + pId)=\rho {\vec g}\cr
\mbox{div}({\vec u}+L_0V_a\nabla Z(\rho))=0.\ea\es\ee
A  stationnary laminar solution satisfies
$$\bs\ba{l}\rho_0(x)u_0(x)=-\rho_aV_a\cr
{d\over dx}(\rho_0(x)u_0(x)^2+p_0(x))=-\rho_0(x)g\cr
{d\over dx}(u_0(x)-L_0V_aZ(\rho_0(x)))=0.\ea\es$$
where the first relation is a consequence of the mass conservation equation, the constant $\rho_0u_0$ being identified through its limit at $x\rightarrow+\infty$. Hence $u_0(x) = -{\rho_aV_a\over \rho_0(x)}$ leading to the equation on $\rho_0$:

$$-{V_a\rho_a\over \rho_0(x)}+L_0V_a{d\rho_0(x)\over dx}\rho_0(x)^{-\nu -2}\rho_a^{\nu+1}=C_0.$$
As $y={x\over L_0}$ this equation rewrites
$${d\xi \over dy}= {C_0\over V_a}\xi^{\nu+2} + \xi^{\nu+1}.$$
If $C_0=0$, the equation becomes ${d\over dy}(\xi^{-\nu})=-\nu$, hence $\xi^{-\nu}=D_0-\nu y$ hence $\xi$ is not defined for $y>\frac{D_0}{\nu}$. We cannot consider this solution.\\
If $C_0>0$, $\xi$ is increasing, hence if it is majorated, it has a limit $l>0$ when $y\rightarrow +\infty$, this limit $l$ satisfies $l^{\nu+1}+{C_0\over V_a}l^{\nu+2}=0$, which is impossible.
We deduce that $C_0<0$, hence the equation is
$${d\xi\over dy}=\xi^{\nu+1}(1-{|C_0|\over V_a}\xi)$$
hence from the resolution of the equation we deduce $\xi\rightarrow {V_a\over |C_0|}$, and as $\xi\rightarrow 1$, $|C_0|=V_a$ and the resulting equation is (\ref{xi}).

This equation has a unique constant solution $\xi=1$. The low Mach approximation of S. Dellacherie \cite{Della} for a bubble model uses this stationnary solution as base solution. When we consider a non constant solution, we have the following Lemma which gives the behavior of the solution when $x$ goes to $-\infty$. Introduce the function $h_{\{\nu\}}$ such that
 $h_{\{\nu(\xi)\}}=\int_0^{\xi}{1\over \eta^{\nu -[\nu]}(1-\eta)}$, and let  $n=[\nu]$, 
${\nu+1\over \nu+2}=\xi_*$ and $y_*=\sum_{p=0}^n{1\over \nu-p}\xi_*^{\nu-p} - h_{\{\nu\}}(\xi_*)$.

\begin{lemma}
There exists $t_*>0$ such that for $0\leq t\leq t_*$ there exists a unique continuous function $g(t)$, $g(0)=1$ solution of

$$(g(t))^{\nu}=-t^{\nu}g(t)^{\nu}(h_{\nu}(tg(t))+y_0)+\sum_{p=0}^nt^p(g(t))^p{\nu\over \nu - p}.$$
The function $g$ has a Taylor expansion at $t=0$.
\begin{enumerate}
\item Behavior when $y\rightarrow -\infty$\\
For $y\leq -{\nu\over (t_*)^{\nu}}$, we have the identity

$$\xi(y)=(-{1\over \nu y})^{1\over \nu}g((-{1\over \nu y})^{1\over \nu}).$$
There exists a function $r(t, \varepsilon)$ such that
\be\label{defifonctionr}\frac{1}{\varepsilon}(\xi(-\frac{t}{\varepsilon}))^{\nu}= \frac{1}{\nu t}+\varepsilon^{\frac{1}{\nu}}t^{-1-\frac{1}{\nu}}r(t, \varepsilon).\ee
There exists $t_0>0$ and $\varepsilon_0>0$ such that $r(t, \varepsilon)$ is bounded for $t\geq t_0, 0\leq\varepsilon\leq \varepsilon_0$
\item Auxiliary function at $-\infty$:
We introduce $S(t, \varepsilon)=-\int_t^{+\infty}s^{-1-\frac{1}{\nu}}r(s, \varepsilon)ds$. This function satisfies $S(t, \varepsilon)t^{\frac{1}{\nu}}$ uniformy bounded for $t\geq t_0$ and $0\leq \varepsilon\leq \varepsilon_0$.
\item Behavior when $y\rightarrow +\infty$\\
We have $1-\xi(y)=C(y)e^{-y}$, with $C(y)\rightarrow exp(y_0-{1\over \nu}-..-{1\over \nu-n}+\int_0^1{(1-\eta^{\nu-n})d\eta\over \eta^{\nu-n}(1-\eta)})$.
\end{enumerate}
\label{comportementxi}
\end{lemma}
The proof of this Lemma uses the following relation

$$-\sum_{p=0}^n{1\over \nu - p}\xi^{-\nu + p} + h_{\{\nu\}}(\xi)=y-y_*$$
The equality yields
$$y=-{1\over \nu \xi^{\nu}}[\sum_{p=0}^n{\nu\over \nu-p}\xi^p-\nu\xi^{\nu}(h_{\{\nu\}}(\xi)+y_*)].$$
Introduce $\xi=tg$ and $t=(-{1\over \nu y})^{1\over \nu}$, such that $y=-{1\over \nu t^\nu}$. We obtain the equality

$$g^{\nu}=\sum_{p=0}^n{\nu\over \nu - p}t^pg^p-\nu t^{\nu}g^{\nu}[h_{\{\nu\}}(tg)+y_*].$$
We introduce $\Psi(t, g)=g^{\nu}-\sum_{p=0}^n{\nu\over \nu-p}t^pg^p-\nu t^\nu g^\nu[h_{\{\nu\}}(tg)+y_0]$. We have $\Psi(0, 1)=0$ and $\partial_g\Psi(0, 1)=\nu>0$ hence by the implicit function theorem, there exists a unique solution of $\Psi(t, g)=0$ in the neighborhood of $(1, 0)$.

The behavior of the solution is more classical in the neighborhood of $+\infty$, thanks to the equality

$$\ln(1-\xi)=y_0-\sum_{p=0}^n{1\over (\nu - p)\xi^{\nu-p}} +\int_0^{\xi}{d\eta (1-\eta^{\nu - n})\over \eta^{\nu -n}(1-\eta)} -y$$
hence
$$1-\xi(y)=e^{C(y)}e^{-y}$$
with $C(y)\rightarrow y_0-\sum_{p=0}^n{1\over (\nu - p)} +\int_0^1{d\eta (1-\eta^{\nu - n})\over \eta^{\nu -n}(1-\eta)}$ when $y\rightarrow +\infty$. The Lemma is proven.
\paragraph{The particular case $\nu= 2.5$}
In this case, explicit calculations lead to the following implicit relation

$$y-C=\ln\frac{1+\sqrt\xi}{1-\sqrt\xi}-\frac{2}{5\xi^\frac52}(1+\frac53 \xi+5\xi^2).$$
Hence, for $y=-\frac{t}{\varepsilon}$ and $\xi=\eta \varepsilon^{\frac{2}{5}}$ we get

$$t=\frac{2}{5\eta^{\frac{5}{2}}}(1+\varepsilon^{\frac{2}{5}}\eta+ \varepsilon^{\frac{4}{5}}\eta^2)-\varepsilon C-\varepsilon\ln\frac{1+\varepsilon^{\frac{1}{5}}\sqrt\eta}{1-\varepsilon^{\frac{1}{5}}\sqrt\eta}.$$
We write $\eta(\varepsilon, t)=(\frac{5t}{2})^{-\frac{2}{5}}g(t, \varepsilon)$. Considering the limit, for $t>0$ fixed, in the previous equality we obtain
$$\mbox{lim}_{\varepsilon\rightarrow 0}g(t, \varepsilon)=1$$
and we construct step by step the expansion of $g(t, \varepsilon)$ in $\varepsilon$.

\subsection{Physical interpretation of the model}
We insist finally on the fact that this system of equations is only a theoretical model: the equation satisfied by $p_0$ is
$${dp_0\over dy}=-\rho_aV_a^2[{1\over Fr}\xi(y)+\xi^{\nu-1}(1-\xi)]$$
hence the pressure is not bounded. The total pressure writes

$$P(x)=(C_p-C_v)\rho_a T_a + M^2p_0(L_0x)$$
which leading order term is constant and low order term in Mach is not bounded. Hence a good way of calling this model could be to call it a {\bf relative isobaric model or a relative low mach model. } See Majda and Sethian \cite{Majda}, Embid \cite{ME}, or P.L. Lions \cite{lions} for other remarks on this modelling.\\
Introduce the function $M(y)$ such that $M^2(y)={|{\vec u}(y)|^2\over C_s^2(y)}$, where $C_s(y)$ is the sound velocity at a point of the fluid given by $C_s(y)^2=\frac{C_p}{C_v}\frac{p(y)}{\rho(y)}$, that we still call the Mach number. We have
\begin{lemma}
The Mach number of the stationnary solution is bounded for $\nu\geq 2$ and $y\leq -C$.
\end{lemma}
\paragraph{Proof} As $M^2(y)={\rho_a^2 V_a^2\over \xi(y)C_s^2(y)}={\rho_aV_a^2\over \gamma p(y)\xi(y)}= {\rho_aV_a^2\over (C_p-C_v)\rho_aT_a\xi(y) + M^2p_0(y)\xi(y)M^2}$
which has a finite limit when $y\rightarrow -\infty$ under the condition $\nu>2$. Moreover, as $p_0(y)\rightarrow -\infty$ when $y\rightarrow +\infty$, there exists a point where the pressure vanishes. There exists a constant $C$ such that $P(y)$ is bounded below $\varepsilon_0>0$ on $]-\infty, -C]$ and $-C$ is a positive constant of order $({\rho_aT_a\over M^2})^{\nu\over \nu-1}$.

The low Mach number assumption is relevant (in particular when looking at the temperature of ablation and the density in the ablated fluid).
\subsection{Linearisation of the equations}
To simplify the notations of what follows, we denote by ${\tilde f}$ a quantity appearing in the Euler system of equations, by $f_0$ its stationnary leading order term, and by $f$ the perturbation of order 1 normalized by the physical quantity $\rho_a$ for the density, $\rho_aV_a$ for the impulsion, and $\rho_aV_a^2$ for a pressure term.
In the normal mode study in the vicinity of a profile depending on $x$, it is pertinent to linearize the variables on which acts the derivative ${d\over dx}$, that is ${\tilde \rho}{\tilde u}$, ${\tilde \rho}({\tilde u})^2+p$ and ${\tilde \rho}{\tilde u}{\tilde v}$, which writes
\be\label{linearized1}\bs\ba{l}{\tilde \rho}{\tilde u}=-\rho_aV_a+\rho_aV_ax_1\cr
{\tilde \rho}({\tilde u})^2+p=p_0(x)+\rho_0(x)(u_0(x))^2+\rho_aV_a^2x_2\cr
{\tilde \rho}{\tilde u}{\tilde v}=-i\rho_aV_a^2x_3\ea\es\ee

The coupling with the equations of the energy is made through the perturbation of density $\rho$. Introduce
$${\vec \tau}={\vec {\tilde u}}-L_0V_a\nabla (Z(\rho_0)) - L_0V_a\nabla (Z({\tilde \rho})-Z(\rho_0)).$$
The thermal perturbed quantities are \be \label{linearized2}x_4=Z(\rho_0) -Z({\tilde \rho}), x_5={1\over V_a}(\tau_1+V_a).\ee
With the choice of unknowns $(x_1, x_2, x_3, x_4, x_5)$ (which corresponds to the unknowns $\rho u, \rho u^2+p, \rho u v$, $-Z(\rho), \tau_1$), denoting by $Z^{-1}$ the inverse function of $Z$, $Z^{-1}(f)=\rho_a((\nu+1)f)^{1\over \nu+1}$, the non linear system (\ref{KANL}) is equivalent to:
\be\bs\ba{l}
V_a^{-1}\partial_t(Z^{-1}(Z(\rho_0)-x_4))+\rho_a\partial_xx_1+i\partial_z(Z^{-1}(Z(\rho_0)-x_4){x_3\over 1-x_1})=0\cr
V_a^{-1}\partial_tx_1+\rho_a\partial_xx_2-i\rho_a\partial_z x_3 = -{g\over V_a^2}(Z^{-1}(Z(\rho_0)-x_4))+{g\over V_a^2}\rho_0(x)\cr
V_a^{-1}\partial_t(Z^{-1}(Z(\rho_0)-x_4){x_3\over 1-x_1})-\rho_a\partial_xx_3\cr
+\rho_a\partial_z(x_2+{x_3^2\over (1-x_1)^2}{Z^{-1}(Z(\rho_0)-x_4)\over \rho_a}+{\rho_a\over \rho_0(x)}-{\rho_a(-x_1)^2\over Z^{-1}(Z(\rho_0)-x_4)})=0\cr
L_0\partial_xx_4 +{\rho_a(1-x_1)\over Z^{-1}(Z(\rho_0)-x_4)}-{\rho_a\over \rho_0(x)}-x_5=0\cr
\partial_xx_5+\partial_z(L_0\partial_zx_4+i{x_3\over 1-x_1})=0.\ea\es\ee
In the system (\ref{KANL}), one needs to obtain the linearization of ${\tilde \rho}{\tilde v}$ and of ${\tilde p}+{\tilde \rho}({\tilde v})^2$ in terms of $x_1, x_2, x_3$. It is a consequence of
$${\tilde v}=iV_a{x_3\over 1-x_1}, {\tilde \rho}=\rho_a (\xi+\rho), {\tilde p}=x_2-{(1-x_1)^2\over \xi+\rho}+{1\over \xi}, {\tilde u}=-\frac{V_a}{\xi+\rho)}$$
hence the approximations (dropping the terms of order 2 at least)
$$ {\tilde v}=iV_ax_3, {\tilde p}=p_0 + \rho_aV_a^2(x_2+{2\over \xi}x_1+{\rho\over \xi^2}), {\tilde u}=u_0+\frac{V_a}{\xi^2}\rho.$$
Moreover, the identity $x_4=Z(\rho_a\xi)-Z(\rho_a(\xi+\rho))$ leads to $x_4=-Z'(\rho_0)\rho_a\rho+O(\rho^2)$, hence
$$x_4={\rho\over \xi^{\nu+2}}+O({\rho^2\over \xi^{\nu+4}}).$$
Moreover, as
$$\tau_1={\tilde u}-u_0-V_a-L_0V_a\partial_x(Z(\rho_a(\xi+\rho))-Z(\rho_a\xi))$$
one deduces
$$x_5=V_a^{-1}({\tilde u}-u_0)-\rho_aL_0\partial_x(Z'(\rho_a\xi)\rho).$$
The linearized system is
\be\bs\ba{l}V_a^{-1}\partial_t\rho+\partial_xx_1+\xi\partial_z(ix_3)=0\cr
V_a^{-1}\partial_tx_1+\partial_xx_2-\partial_z(ix_3)=\frac{g}{V_a^2}\rho\cr
V_a^{-1}\partial_t(i\xi x_3)-\partial_x(ix_3)+\xi^{-1}\partial_z(2x_1-\xi^{-1}\rho)=0\cr
\rho= -(\rho_aZ'(\rho_a\xi))^{-1}x_4\cr
L_0\partial_x(\rho_aZ'(\rho_a\xi)\rho)+x_5+ \frac{V_a}{\xi^2}\rho=0\cr
\partial_xx_5+\partial_z(L_0\partial_zx_4+ix_3)=0\ea\es\label{lspde}\ee
Using the new variable $y$ such that $x=L_0y$, one deduces the relation
$$ x_5=\partial_yx_4+{\rho + \xi(y)x_1\over \xi(y)(\xi(y)+\rho)}= \partial_yx_4+\xi^\nu x_4+{x_1\over \xi(y)} + O({x_1\rho + \rho^2\over \xi^2}).$$
Write the following normal mode expression:
\be\label{normal}\left(\ba{c}x_1\cr
x_2\cr
ix_3\cr
x_4\cr
x_5\ea\right)=\Re (\left(\ba{c}{\underline x}_1\cr
{\underline x}_2\cr
i{\underline x}_3\cr
{\underline x}_4\cr
{\underline x}_5\ea\right)e^{ikz+\gamma\sqrt{gk}t}).\ee
Assume $\left(\ba{c}{\underline x}_1\cr
{\underline x}_2\cr
i{\underline x}_3\cr
{\underline x}_4\cr
{\underline x}_5\ea\right)e^{ikz+\gamma\sqrt{gk}t}$ is a solution of the linearized system. This rewrites as a system of ordinary differential equations  on$\left(\ba{c}{\underline x}_1\cr
{\underline x}_2\cr
i{\underline x}_3\cr
{\underline x}_4\cr
{\underline x}_5\ea\right)$. If $\left(\ba{c}{\underline x}_1\cr
{\underline x}_2\cr
i{\underline x}_3\cr
{\underline x}_4\cr
{\underline x}_5\ea\right)$ is solution of this linear system of ODE, then $(x_1, x_2, ix_3, x_4, x_5)^{t}$ is solution of (\ref{lspde}). Note that, in this case $ix_3$ is real.\\
\paragraph{Remark on complex growth rates}
Note that, if $\gamma$ is complex, the solution $\left(\ba{c}{\underline x}_1\cr
{\underline x}_2\cr
i{\underline x}_3\cr
{\underline x}_4\cr
{\underline x}_5\ea\right)$ is also complex and depend on $y, k, \gamma$. More precisely, introduce $y_j, z_j$ such that
${\underline x_j}=y_j+iz_j$. We have thus
$$\Re (\left(\ba{c}{\underline x}_1\cr
{\underline x}_2\cr
i{\underline x}_3\cr
{\underline x}_4\cr
{\underline x}_5\ea\right)e^{ikz+\gamma\sqrt{gk}t})\vert_{t=0}=(\left(\ba{c}y_1\cos kz-z_1\sin kz\cr
y_2\cos kz -z_2\sin kz\cr
-z_3\cos kz -y_3\sin kz\cr
y_4\cos kz-z_4\sin kz\cr
y_5\cos kz-z_5\sin kz\ea\right)$$
and 
$y_j$ and $z_j$ are known through the decomposition of the initial perturbation in the even part and the odd part in $z$.\\
Hence, from a complex solution $({\underline x_1}, ..., {\underline x_5})^{t}(x, k, \gamma)$ of the normal mode system, one deduces a solution of the perturbation system with a known initial condition.\\
\paragraph{Normal modes system}
The equations for the normal modes associated with the mass conservation and the momentum equation are
\be\bs\ba{l}\alpha\gamma\rho + {dx_1\over dy}+\alpha\beta\xi x_3=0\cr
{dx_2\over dy}-\alpha\beta x_3+\alpha\gamma x_1+{\alpha\over \beta}\rho=0\cr
{dx_3\over dy}-\alpha\gamma\xi x_3+\alpha\beta (x_2+{2\over \xi}x_1+{\rho\over \xi^2})=0.\ea\es\ee
The normal mode formulation of the linearized energy equation is
$${dx_5\over dy}+i\alpha\beta(-ix_3+i\alpha\beta x_4)=0.$$
Let $X$ be given by
$$X=\left(\ba{l}x_1\cr x_2\cr x_3\cr x_4\cr x_5\ea\right).$$
The linearized system on $X$ is
\be \label{systemefonda}{dX\over dy}+M_0(\xi, \alpha,\beta, \gamma)X=0\ee
where the matrix $M_0$ is given by
\be M_0(\xi, \alpha, \beta, \gamma)=\left(\ba{ccccc}0&0&\alpha\beta\xi&\alpha\gamma\xi^{\nu+2}&0\cr
\alpha\gamma&0&-\alpha\beta&{\alpha\over \beta}\xi^{\nu+2}&0\cr
{2\alpha\beta\over \xi}&\alpha\beta&-\alpha\gamma\xi&\alpha\beta\xi^{\nu}&0\cr
{1\over \xi}&0&0&\xi^{\nu}&-1\cr
0&0&\alpha\beta&-\alpha^2\beta^2&0\ea\right).\ee
{\bf From now on, we will call this system the Kull-Anisimov
  system.}
The eigenvalues of $-M_0$ are given by the classical result\footnote{A
  change of unknowns (that is a general ${\tilde X}=R(y)X$) lead to a
  different set of eigenvalues, however if the matrix $R(y)$ depend
  only on $\xi(y)$ and is a $C^1$ function of $\xi\in [0, 1]$ then the
  limit of the eigenvalues when $y\rightarrow \pm \infty$ is the same
  as the limit of the eigenvalues described above.} (see \cite{Kull},
\cite{Piriz}), and help us to study the solution at $\pm \infty$:
\begin{prop}
\begin{itemize}
\item The eigenvalues of $-M_0$ are
\be\label{valpropres}\lambda_0(\xi)=\alpha\gamma\xi, \lambda_{a,
  +}=\alpha\beta, \lambda_{a, -}=-\alpha\beta, \lambda_+(\xi),
\lambda_-(\xi)\ee
where (the square root is chosen of positive real part)
\be\lambda_{\pm}(\xi)= -{\xi^{\nu}\over 2}\pm \sqrt{{\xi^{2\nu}\over 4}+\alpha\gamma\xi^{\nu+1}+\alpha^2\beta^2}\ee

The eigenvalues $\lambda_0, \lambda_{a, \pm}$ are called the
hydrodynamic modes, the eigenvalues $\lambda_\pm(\xi)$ are called the
thermal modes.\\
\item For $\Re\gamma\geq 0$, one has $\Re \lambda_0(\xi)\geq 0$,
  $\pm\Re(\lambda_\pm(\xi))>0$, and for $\Re\gamma>0$ one has $\Re
  \lambda_0(1)>0$. The matrix $-M_0(\xi)$ has three eigenvalues of
  positive real part, and two eigenvalues of negative real part. 
\item The associated eigenvectors are given by

\be\label{vecteurspropreshydro}\ba{l}E_0(\xi)=\beta i -2\beta
i_2-\gamma\xi i_3\cr
E_{a, +}(\xi)=-\beta i + (\beta +\gamma\xi)i_2+\beta i_3\cr
E_{a, -}(\xi)=-\beta i + (\beta -\gamma\xi)i_2-\beta i_3\ea\ee
If we introduce
$${\tilde R}^0=\left(\ba{c}\gamma\xi^2\cr {\xi^2\over \beta}\cr \beta\cr 0\cr 0\ea\right), {\tilde T}^0=\left(\ba{c}\xi(\gamma^2\xi^2-\beta^2)\cr \beta^2-\gamma^2\xi^2+{\gamma\xi^3\over \beta}\cr -\xi^2\cr 0\cr 0\ea\right)$$
the eigenvectors $F_{\pm}$ associated with $\lambda_{\pm}(\xi)$ are given by

\be\label{f}F_{\pm}=i_4+\alpha{\tilde R}^0 + {\alpha^2\over
  \lambda_\pm-\alpha\gamma\xi}{\tilde T}^0+{\alpha^3\xi^2\over
  (\lambda_\pm-\alpha\gamma\xi)^2}E_0(\xi).\ee
\end{itemize}
\end{prop}
The proof of this Proposition is straightforward, except for the sign
of the real part of $\lambda_\pm(\xi)$. For this, we use
$$(\lambda+\frac{\xi^{\nu}}{2})^2=\frac{\xi^{2\nu}}{4}+\alpha^2\beta^2+\alpha\xi(\Re\gamma
+i\Im \gamma)=(A+iB)^2$$
where $A>0$ (if $A<0$ one uses $(A+iB)^2=(-A-iB)^2$, and $A=0$ leading
to $(iB)^2=-B^2$ is not possible when $\Re \gamma\geq 0$).
Hence one obtains
$$A^2-B^2=\frac{\xi^{2\nu}}{4}+\alpha^2\beta^2+\alpha\xi\Re\gamma,
A^2+B^2=((\frac{\xi^{2\nu}}{4}+\alpha^2\beta^2+\alpha\xi\Re\gamma)^2+\alpha^2\xi^2(\Im
\gamma)^2)^{\frac12}$$
As $\lambda_\pm(\xi)=-\frac{\xi^{\nu}}{2}\pm(A+iB)$, $\Re
(\lambda_+(\xi)\lambda_-(\xi))=(-\frac{\xi^{\nu}}{2}+A)(-\frac{\xi^{\nu}}{2}-A)+B^2$,
and
$\Re\lambda_+(\xi)\Re\lambda_-(\xi)=\Re(\lambda_+(\xi)\lambda_-(\xi))-B^2=-\alpha\xi^{\nu+1}\Re\gamma
-\alpha^2\beta^2-B^2<0$, one obtains that the product
$\Re\lambda_+(\xi)\Re\lambda_-(\xi)$ is strictly negative, hence the real parts
are of opposite sign, hence $\Re\lambda_+(\xi)=-\frac{\xi^{\nu}}{2}+A$.\\
This calculus also defines uniquely in the case $\Re\gamma>0$ the
eigenvalues $\lambda_+(\xi)$ and $\lambda_-(\xi)$.\\
Note also $A_0(\alpha, \beta, \gamma)$ and $B_0(\alpha, \beta, \gamma)$ the quantities such that $A_0(\alpha, \beta, \gamma)>0$ and $(A_0(\alpha, \beta, \gamma)+iB_0(\alpha, \beta, \gamma))^2=\frac14+\alpha^2\beta^2+\alpha\gamma$.\\
Using ${d\xi\over dy}=\xi^{\nu+1}(1-\xi)$, we also introduce a new system of unknowns :
\be Y = TX=\left(\ba{c}x_1-\alpha\gamma\xi x_4\cr
x_2-{\alpha\over \beta}\xi x_4\cr
x_3-\alpha\beta x_4\cr
{\xi\over 1-\xi}\alpha\beta x_4\cr
x_5\ea\right)\ee
and the system on $Y$, equivalent to (\ref{systemefonda}) is
$$\bs\ba{l}\frac{dy_1}{dy}-\alpha\gamma y_1+\alpha\beta\xi y_3+\alpha(1-\xi)(\beta^2-\frac{\gamma^2}{\beta})z_4+\alpha\gamma\xi x_5=0\cr
\frac{dy_2}{dy}+\alpha(\gamma-\frac{1}{\beta}) y_1-\alpha\beta y_3+\alpha(1-\xi)(\frac{\gamma^2}{\beta}-\frac{\gamma}{\beta^2}-\frac{\beta}{\xi})z_4+\alpha\frac{\xi}{\beta} x_5=0\cr
\frac{dy_3}{dy}+\frac{\alpha\beta}{\xi}y_1+\alpha\beta y_2-\alpha\gamma\xi y_3+\alpha(1-\xi)(\frac{1}{\beta}\frac{\gamma(1-\xi)}{\xi})z_4+\alpha\beta x_5=0\cr
\frac{dz_4}{dy}+\alpha\gamma z_4\frac{\alpha\beta}{1-\xi}(y_1-\xi x_5)=0\cr
\frac{dx_5}{dy}+\alpha\beta y_3=0,\ea\es$$
namely
\be\label{SL2}{dY\over dy}+\alpha B(\xi, \beta, \gamma)Y=0\ee
where
$$\bs\ba{l}Bi_1=-\gamma i_1+(\gamma - {1\over \beta})i_2+{\beta\over \xi}i_3 + {\beta\over 1-\xi}i_4\cr
Bi_2=\beta i_3\cr
Bi_3 = \beta\xi i_1 - \beta i_2 - \gamma\xi i_3 + \beta i_5\cr
Bi_4=(1-\xi)[(\beta - {\gamma^2\over \beta})i_1+({\gamma\over \beta}(\gamma -{1\over \beta})-{\beta\over \xi})i_2+({1\over \beta}+\gamma{1-\xi\over \xi})i_3] + \gamma i_4\cr
Bi_5=\gamma\xi i_1 + {\xi\over \beta} i_2 + \beta i_3 -{\beta\xi\over 1-\xi}i_4.\ea\es$$

The next section is devoted to the statement of the methods used to find the solutions at infinity for systems of ODE which coefficients depend on $y$, and of the general set-up to find solutions bounded at $\pm \infty$.
\section{Evans functions and application to the Kull-Anisimov system}
\label{evan}
\subsection{General framework}
This section recalls results of the paper of Alexander, Gardner
et Jones \cite{agj}, as well as the methods developed by K. Zumbrun \cite{Zum},
D. Serre \cite{Zum1}, S. Benzoni-Gavage \cite{sbg} and other authors.  Its purpose is to study solutions with a prescribed behavior at infinity of an ordinary linear system of differential equations. It is used in particular to identify solutions going to 0 as $y\rightarrow \pm \infty$.\\

In the general case, we consider the ordinary differential system
\be\label{Esystem}{dy\over dt}=A(t, \alpha)y,\ee
when $A$ is a regular matrix (for example analytic in $\alpha$).\\
We notice that the vectorial product $y_1\wedge y_2$ of two solutions $y_1$ and $y_2$ of (\ref{Esystem}) is solution of a new differential system on $\Lambda^2(K^n)$ which matrix is denoted by $A^{(2)}$, because
$${d\over dt}(y_1\wedge y_2)=(Ay_1\wedge y_2+y_1\wedge Ay_2).$$
Similarily $y_1\wedge y_2\wedge... \wedge y_k$ is solution of an ordinary differential system on $\Lambda^k(\RR^n)$ whose matrix is denoted by $A^{(k)}$.
The matrix $A^{(k)}$ is given by
\be\label{vecto}A^{(k)}(e_{i_1}\wedge e_{i_2}...\wedge e_{i_k})=\sum_{l=1}^k e_{i_1}\wedge .. \wedge Ae_{i_l}\wedge .. e_{i_k}.\ee
When the matrix $A$ is diagonalizable, with eigenvalues $\lambda_1\leq\lambda_2<...\leq\lambda_d$ then $A^{(k)}$ is diagonalizable and its eigenvalues are

$$\sum_{i\in I, I\subset\{1,..., d\}, \mbox{Card}(I)=k}\lambda_i.$$
The largest eigenvalue of $A^{(k)}$ is $\sum_{p=1}^k\lambda_{d+1-p}$, and its smallest eigenvalue is $\sum_{p=1}^k \lambda_p$.\\
Under the hypothesis that $\lambda_{d-k}<\lambda_{d+1-k}$, the largest eigenvalue of $A^{(k)}$ is simple. Its associated eigenvector is the vectorial product of the eigenvectors associated with $(e_{d+1-k},.., e_d)$.

Recall that the space $\Lambda^{(d)}(\RR^d)$ is of dimension 1, hence the matrix associated with $A^{(d)}$ is a number, which is equal to $\mbox{Tr}(A(t, \alpha))$. The associated differential equation is
$${d\over dt}(y_1\wedge ...\wedge y_d)=\mbox{Tr}(A(t, \alpha))(y_1\wedge ...\wedge y_d).$$
Hence the vectorial product of  $d$ solutions of (\ref{Esystem}) satisfies
\be\label{32bis}y_1\wedge y_2\wedge .. \wedge y_d(t, \alpha)=
y_1\wedge y_2\wedge .. \wedge y_d(t_0, \alpha)exp(\int_{t_0}^t\mbox{Tr}
(A(s, \alpha))ds).\ee
whose solution is the Wronskian of $d$ solutions of the system.
\subsection{Notations for the Kull-Anisimov system}
{\bf Remark that, when $\gamma$ is complex, one has to replace $\RR$ by $K=\CC$ but nothing will change as what is important is that we study objects on a field, which can be $\RR$ when $\gamma$ is real, and which is $\CC$ when $\gamma$ is complex.}

In the set-up of this paper, the matrix $-M_0$ admits three eigenvalues of positive real part, which may be associated with the solutions going to 0 when $y$ goes to $+\infty$, and has two eigenvalues of negative real part, which help to understand the solutions going to 0 when $y$ goes to $+\infty$.

It is in general hard to compute the solutions associated with an eigenvalue of the matrix $-M_0$. However, we may compute the solution for the matrix $M_0^{(2)}$ associated with the smallest eigenvalue $\lambda_{a, -}+\lambda_-(\xi)$, and the solution for the matrix $M_0^{(3)}$ associated with the largest eigenvalue $\lambda_0(\xi)+\lambda_+(\xi)+\lambda_{a, +}$.\\
We introduce from now on the base vectors in $\Lambda^2(K^5)$ which take into account the role of ${1\over 1-\xi}$ in $B$:

\be\ba{c}f_1=i_1\wedge i_4, f_2=i_2\wedge i_4, f_3=i_3\wedge i_4, f_4=i_4\wedge i_5\cr
\ba{l}g_1=i_1\wedge i_2, g_2=i_1\wedge i_3, g_3=i_1\wedge i_5\cr
g_4=i_2\wedge i_3, g_5=i_2\wedge i_5, g_6=i_3\wedge i_5\ea\ea\label{vecteursLambda2}\ee
To these vectors are associated the following vectors in $\Lambda^3(K^5)$ such that $f_i\wedge f_i^{\perp}=i_1\wedge i_2\wedge i_3\wedge i_4\wedge i_5=g_j\wedge g_j^{\perp}$, $\forall i, j$. We have

\be\ba{c}f_1^{\perp}=i_2\wedge i_3\wedge i_5, f_2^{\perp}=-i_1\wedge i_3\wedge i_5, f_3^{\perp}=i_1\wedge i_2\wedge i_5, f_4^{\perp}=i_1\wedge i_2\wedge i_3\cr
\ba{l}g_1^{\perp}=i_3\wedge i_4\wedge i_5, g_2^{\perp}=-i_2\wedge i_4\wedge i_5, g_3^{\perp}=-i_2\wedge i_3\wedge i_4\cr
g_4^{\perp}=i_1\wedge i_4\wedge i_5, g_5^{\perp}=i_1\wedge i_3\wedge i_4, g_6^{\perp}=-i_1\wedge i_2\wedge i_4.\ea\ea\label{vecteursLambda3}\ee
Note that we shall use in the sequel the eigenvector of the matrix $M_0^{(2)}(1)$ associated with the eigenvalue $-\alpha \beta + \lambda_-(1)$:
\be\label{vecpro}W_+=F_-(1)\wedge E_{a, -}(1)=\beta f_1+(\gamma - \beta)f_2+\beta f_3-\beta f_4 + G\ee
where $G$ belongs to the space generated by $g_j, j=1..6$, and $W_{+, 1}=\beta$, $W_{+, 2}=\gamma - \beta$, $W_{+, 3}=\beta$, $W_{+, 4}=-\beta$. We also introduce $\mu(\alpha)=-\beta +\frac{\gamma + \alpha\beta^2}{\lambda_-(1)}$ such that
\be\label{defimu} \lambda_{-}(1)-\alpha\beta = -1+\alpha\mu(\alpha).\ee
We notice that $\mu(\alpha)=-\beta -\frac{\gamma+\alpha\beta^2}{\frac12+A_0+iB_0}$, hence \be\label{ine0mu}\forall \alpha\in [0, \alpha_0], \beta\in [\beta_0, \beta_0^{-1}], |\gamma|\leq \beta_0^{-1}, \Re\gamma\geq 0, |\mu(\alpha)|\leq 2\alpha_0\beta_0^{-2}+3\beta_0^{-1}.\ee

We construct a solution of the system
\be\label{s2}{dX^{(2)}\over dy}+M_0^{(2)}X^{(2)}=0\ee
which belongs to the family of its most decreasing solutions when $y$ goes to 
$+\infty$. Similarily, we construct a solution of the system

\be\label{s3}{dX^{(3)}\over dy}+M_0^{(3)}X^{(3)}=0\ee
belonging to the family of its most decreasing solution when $y$ goes to 
$-\infty$.\\
It is useful to introduce the following transformation for the study of the solution when $y\rightarrow +\infty$: \be Y^{(2)}= T^{(2)}X^{(2)}\ee
where
\begin{lemma} \label{transfo}
Let $\sum x_ji_j$ and $\sum t_ji_j$ be two solutions of (\ref{systemefonda}). Ve denote by $X^{(2)}=v_1f_1+v_2f_2+v_3f_3+v_4f_4+\sum_{j=1}^6w_jg_j$ a solution of (\ref{SL2}). The associated solution $Y^{(2)}=T^{(2)}X^{(2)}$ writes
$$Y^{(2)}= Z_1f_1+Z_2f_2+Z_3f_3+Z_4f_4 + \sum_{j=1}^6 M_jg_j$$
with
$$\bs\ba{l}Z_j=\alpha\beta {\xi\over 1-\xi}v_j\cr
M_1={\alpha\over \beta}\xi v_1-\alpha\gamma\xi v_2 + w_1\cr
M_2=\alpha\beta v_1-\alpha\gamma\xi v_3+w_2\cr
M_3=\alpha\gamma\xi v_4 + w_3\cr
M_4=\alpha\beta v_2 -{\alpha\over \beta}\xi v_3+w_4\cr
M_5={\alpha\over \beta}\xi v_4+w_5\cr
M_6=\alpha\beta v_4+ w_6.\ea\es$$
\end{lemma}
For the construction of the solutions when $y\rightarrow -\infty$,
we use the following transformation of the unknowns
\be\label{nouvellesinc}\xi z_1=x_1, z_4=\alpha\beta x_4.\ee The system deduced from
(\ref{systemefonda}) is thus
\be\label{systemfuchscompletinitial}\bs\ba{l}\frac{dz_1}{dy}+\xi^{\nu}(1-\xi)z_1+\alpha\beta
x_3+\frac{\gamma}{\beta}\xi^{\nu+1}z_4=0\cr
\frac{dx_2}{dy}+\alpha\gamma\xi z_1-\alpha\beta
x_3+\frac{1}{\beta^2}\xi^{\nu+2}z_4=0\cr
\frac{dx_3}{dy}+2\alpha\beta z_1+\alpha\beta x_2-\alpha\gamma\xi
x_3+\xi^{\nu}z_4=0\cr \frac{dz_4}{dy}+\alpha\beta
z_1+\xi^{\nu}z_4-\alpha\beta x_5=0\cr \frac{dx_5}{dy}+\alpha\beta
x_3-\alpha\beta z_4=0\ea\es\ee Introduce $$ t=-\alpha\beta y.$$ We have
\begin{lemma} \label{transfo-}
Let $Z^1(t, \alpha), Z^2(t, \alpha), Z^3(t, \alpha)$ be three solutions of
(\ref{systemfuchscompletinitial}). Write \be Z^{(3)}(t,
\alpha)=Z^1(t, \alpha)\wedge Z^2(t, \alpha)\wedge Z^3(t, \alpha)=\sum_{j=1}^4f_j(t, \alpha)f_j^{\perp}+\sum_{p=1}^6g_p(t,
\alpha)g_p^{\perp}.\ee The solution of (\ref{s3}) associated with
$Z^{(3)}$ is
$$\ba{ll}w^{(3)}(y, \alpha)&=f_1(-\alpha\beta y, \alpha)f_1^{\perp}+\xi(y)\sum_{j=2}^4f_j(-\alpha\beta y, \alpha)f_j^{\perp}\cr
&+\frac{1}{\alpha\beta}[\sum_{p=1}^3g_p(-\alpha\beta y, \alpha)g_p^{\perp}+\xi\sum_{p=4}^6g_p(-\alpha\beta y, \alpha)g_p^{\perp}].\ea$$
\end{lemma}
The proof of these two lemmas is straightforward.
\subsection{Construction of the Evans function for the Kull-Anisimov system}
Recall that the Evans function is characterized by the vectorial product of five solutions of the system. As it writes $X^{(2)}\wedge w^{(3)}$, we notice that, from Lemma \ref{transfo} and Lemma \ref{transfo-}, we have
\be\label{pvect}\ba{ll}\alpha\beta \xi X^{(2)}\wedge w^{(3)}=&\xi(M_1g_1+M_2g_2+M_3g_3)+\xi^2(M_4g_4+M_5g_5+M_6g_6)\cr
&+(1-\xi)[Z_1(f_1-\frac{\xi}{\beta^2}g_1-g_2)+\xi Z_2 (f_2+\frac{\gamma}{\beta}g_1-g_4)\cr
&+\xi Z_3(f_3+\frac{\gamma}{\beta}g_2-\frac{\xi}{\beta^2}g_4)+\xi Z_4(f_4-\frac{\gamma}{\beta}g_3-\frac{\xi}{\beta^2}g_5-g_6)].\ea\ee
Let $C_0$ be the limit of $(1-\xi(y))e^y$ when $y\rightarrow +\infty$. To ensure uniqueness for the systems (\ref{s2}) and (\ref{s3}) and to adapt the constant in the system on $Z, M$, we consider the solution $w_+^{(2)}$ of the system (\ref{s2}) such that \be\label{conditionunicite2}w_+^{(2)}(y)e^{(\alpha\beta - \lambda_{-}(1))y}\ee converges to $W_+\frac{C_0}{\alpha\beta}$ when $y\rightarrow +\infty$. Similarily, we consider the solution $w_-^{(3)}$ of the system (\ref{s3}) such that
\be\label{conditionunicite3}w_-^{(3)}(-\frac{t}{\alpha\beta}, \alpha)e^{2t+\frac{\gamma}{\beta}\int_{t_*}^t\xi(-\frac{s}{\alpha\beta})ds}\rightarrow 2f_1^{\perp}+f_2^{\perp}+f_3^{\perp}-g_1^{\perp}-g_2^{\perp}-2g_3^{\perp}-g_5^{\perp}-g_6^{\perp}=S.\ee
Let $B^{(2)}$ be given by
\be\bs\ba{ll}B^{(2)}f_1= (\gamma - {1\over \beta})f_2+ {\beta\over \xi}f_3 &+(1-\xi)[({\gamma^2\over \beta}-{\beta\over \xi}-{\gamma\over \beta^2})g_1\cr
&+({1\over \beta}+\gamma{1-\xi\over \xi})g_2]\cr
B^{(2)}f_2=\gamma f_2+\beta f_3 &+(1-\xi)[({\gamma^2\over \beta}-\beta)g_1\cr
&+({1\over \beta}+\gamma{1-\xi\over \xi})g_4]\cr
B^{(2)}f_3=\beta\xi f_1-\beta f_2+\gamma (1-\xi)f_3-\beta f_4&+(1-\xi)[(-\beta +{\gamma^2\over \beta})g_2\cr
&+({\beta\over \xi}+{\gamma\over \beta^2}-{\gamma^2\over \beta})g_4]\cr
B^{(2)}f_4=-\gamma\xi f_1-{\xi\over \beta}f_2 -\beta f_3+\gamma f_4&+(1-\xi)[(\beta -{\gamma^2\over \beta})g_3\cr
&+(-{\beta\over \xi}-{\gamma\over \beta^2}+{\gamma^2\over \beta})g_5\cr
&+({1\over \beta}+\gamma{1-\xi\over \xi})g_6]\cr
B^{(2)}g_1=-{\beta\over 1-\xi}f_2 &-\gamma g_1 +\beta g_2 -{\beta\over \xi}g_4\cr
B^{(2)}g_2=-{\beta\over 1-\xi}f_3&-\beta g_1-\gamma(1+\xi)g_2+\beta g_3\cr
& + (\gamma -{1\over \beta})g_4\cr
B^{(2)}g_3={\beta\over 1-\xi}(-\xi f_1+f_4)&+{\xi\over \beta}g_1+\beta g_2-\gamma g_3\cr
& + (\gamma - {1\over \beta})g_5+ {\beta\over \xi}g_6\cr
B^{(2)}g_4=&-\beta\xi g_1-\gamma\xi g_4+\beta g_5\cr
B^{(2)}g_5=-{\beta\xi\over 1-\xi}f_2&-\gamma\xi g_1 + \beta g_4 + \beta g_6\cr
B^{(2)}g_6=-{\beta\xi\over 1-\xi}f_3&-\gamma\xi g_2 +\beta\xi g_3 -{\xi\over
\beta}g_4\cr
& -\beta g_5 -\gamma\xi g_6.\ea\es\ee
The system on $Y^{(2)}$ is
\be {dY^{(2)}\over dy}+\alpha B^{(2)}Y^{(2)}=0\ee
The Evans function that we shall use is given by:

\begin{Definition}
Introduce

$$Ev(\alpha, \beta, \gamma)=\alpha \beta w_+^{(2)}(y)\wedge w_-^{(3)}(y)e^{-\int_0^y (\alpha\gamma\xi(y')-(\xi(y'))^{\nu})dy'}.$$
This function, independant of $y$, is called the Evans function of the problem.
\end{Definition}
It is easy to derive the
\begin{prop}
The complex number $\gamma$ is an instability growth rate according to Definition \ref{defigrowth} if and only if
$$Ev(\alpha, \beta, \gamma)=0.$$
\label{propdefigrowth}
\end{prop}
The proof of Proposition \ref{propdefigrowth} is to be found in \cite{agj}.
\paragraph{Reduction of the Evans function}
We deduce from the systems (\ref{s2}) and (\ref{s3}) that $w=w_+^{(2)}\wedge w_-^{(3)}$ is solution of (\ref{32bis}). From $\mbox{Tr}M_0= -\alpha\gamma\xi +\xi^{\nu}= -\alpha\gamma \xi + \frac{\dot \xi}{\xi}+\frac{\dot\xi}{1-\xi}$ we deduce that the derivative of $w(y)e^{-\int_0^y (\alpha\gamma\xi(y')-(\xi(y'))^{\nu})dy'}$ is zero, hence
$$\frac{\xi}{1-\xi}w(y)e^{-\int_0^y \alpha\gamma\xi(y')dy'}\mbox{ constant.}$$\\
Using (\ref{pvect}), we obtain
\be\label{evans}\ba{ll}Ev(\alpha, \beta, \gamma)&=[Z_1(f_1-\frac{\xi}{\beta^2}g_1-g_2)+\xi Z_2 (f_2+\frac{\gamma}{\beta}g_1-g_4)\cr
&+\xi Z_3(f_3+\frac{\gamma}{\beta}g_2-\frac{\xi}{\beta^2}g_4)+\xi Z_4(f_4-\frac{\gamma}{\beta}g_3-\frac{\xi}{\beta^2}g_5-g_6)\cr
&+\xi\frac{M_1g_1+M_2g_2+M_3g_3}{1-\xi}\cr
&+\xi^2\frac{M_4g_4+M_5g_5+M_6g_6}{1-\xi}]e^{-\int_0^y \alpha\gamma\xi(y')dy'}\frac{1-\xi(0)}{\xi(0)}.\ea\ee
\subsection{Statement of the principal tools for the study of the Evans function}
The aim of this paper is to compute the Evans function through the calculation\footnote{Note that the techniques of differential equations allow us only to compute the solution of a differential system with non constant coefficients only when this solution is associated with the largest eigenvalue of the matrix in the neighborhood of $-\infty$ and with the smallest eigenvalue of the matrix in the neighborhood of $+\infty$.} of $w_+^{(2)}$ and $w_-^{(3)}$.

\begin{theo}
\label{latotale}
$\bullet$ For all $\xi_0>0$, there exists $\alpha_0>0, \beta_0>0$ such that, for $\alpha\leq \alpha_0, \beta\in [\beta_0, {1\over \beta_0}], \gamma\in [0, {1\over \beta_0}]$ and $y$ such that $\xi(y)\in [\xi_0, 1]$ there exists a unique solution $w_+^{(2)}$ of (\ref{s2}) satisfying (\ref{conditionunicite2}).\\
$\bullet$ Let
\be\label{+}{\tilde w}=T^{(2)}w_+^{(2)}(y)\exp(-\alpha\mu(\alpha)y)- (W_{+, 1}f_1+W_{+, 2}f_2+W_{+, 3}f_3+W_{+, 4}f_4).\ee
The function ${\tilde w}$ is analytic in $(y, \alpha$) for $\xi(y)\in [\xi_0, 1[$ and for $\alpha\leq \alpha_0$. For all $\alpha_0>0$ and for all $\xi_0\in ]0, 1[$ there exists a constant $C(\xi_0, \alpha_0)$ such that $$\forall y, \xi(y)\in ]\xi_0, 1[, |{\tilde w}(y)|\leq C(\xi_0, \alpha_0)(1-\xi(y)).$$
$\bullet$ There exists $\alpha_1\leq \alpha_0$ and $R>0$, depending only on $\xi_0$, such that for $\zeta_0<{1\over R}$, $w$ admits an analytic extension for $y$ such that $y\in [({\alpha\over \zeta_0})^{1\over \nu}, \xi_0]$.\\
$\bullet$ For $\alpha_0>0$ and $\beta_0>0$ given, and for all $\alpha\leq \alpha_0$, $ \beta\in [\beta_0, {1\over \beta_0}], \gamma\in [0, {1\over \beta_0}]$, there exists $t_0>0$ such that there exists a unique solution $w_-^{(3)}(y, \alpha)$ of (\ref{s3}) for $-\alpha\beta y\in [t_0, +\infty[$ satisfying (\ref{conditionunicite3}). Moreover we have
\be\label{estimationunique3}w_-^{(3)}(-{t\over \alpha\beta}, \alpha)\exp(-\int_0^{-{t\over \alpha\beta}}(\lambda_++\alpha(\beta+\gamma\xi(y')))dy')-W_0(t)e^{2t}t^{-{1\over 2\nu}}= O(\alpha^{1\over \nu})\ee
uniformly for $t\in [t_0, +\infty[$.
\end{theo}
We see in this theorem the division in three regions for the computation of the solutions of (\ref{s2}) and of (\ref{s3}). The first region $\xi(y)\in [\xi_0, 1[$ is the aim of Section \ref{over}, the second region (which extends the result of $[\xi_0, 1[$) is studied in Section \ref{suniform}, and the solution in the neighborhood of $-\infty$ is characterized in Section \ref{hypergeometrique}. In the next paragraph, we describe the systems that we shall use in what follows. One of the main problems is that the problem to solve is {\bf not} a Cauchy problem, but we have informations at $\pm \infty$ for the solution of (\ref{systemefonda}) that we want to study. Let us say a word on this system. As $\xi$ goes to 0 when $y$ goes to $-\infty$, we notice that the matrix $M_0$ is singular when $y$ goes to $-\infty$.\\
Moreover, the term ${\alpha\over \beta}$ prevents us to have a result which is uniform when $\beta$ goes to 0. A possible choice to overcome this difficulty is the choice $\beta\in [\beta_0, \frac{1}{\beta_0}]$.\\
Finally, even if we remove the singularity of $M_0$ at $y\rightarrow -\infty$ by whatever method, another problem is induced by the behavior of $\xi$ when $y$ goes to $-\infty$, because $\xi$ goes to 0 as $|y|^{-{1\over \nu}}$. This will lead to a fuchsian problem (Section \ref{hypergeometrique}).

Note that in the theorem \ref{latotale} we defined unique solutions of the problems (\ref{s2}) and (\ref{s3}) with the prescribed behavior at infinity.\\
The behavior of the solutions induced by the theorem \ref{latotale} lead to another expression of the Evans function of relation (\ref{evans}). \\
Assume that $y<0$ and denote by $t=-\alpha\beta y$. We introduce the functions $L_p(t, \alpha), R_j(t, \alpha)$ given by the equalities
\be\label{RL}\ba{c}g_p(t, \alpha)e^{2t+\frac{\gamma}{\beta}\int_{t_0}^t\xi(-\frac{t'}{\alpha\beta})dt'}t^{-\frac{1}{2\nu}}= L_p(t, \alpha), 1\leq p\leq 6\cr
f_j(t, \alpha)e^{2t+\frac{\gamma}{\beta}\int_{t_0}^t\xi(-\frac{t'}{\alpha\beta})dt'}t^{-\frac{1}{2\nu}}=R_j(t, \alpha), 1\leq j\leq 4.\ea\ee
The corrective factor $e^{2t+\frac{\gamma}{\beta}\int_{t_0}^t\xi(-\frac{t'}{\alpha\beta})dt'}t^{-\frac{1}{2\nu}}$ is induced by the relation (\ref{estimationunique3}).
Similarily, we introduce
\be\label{zm0}z_j(y, \alpha)=Z_j(y, \alpha)e^{-\alpha\mu(\alpha)y}, m_p(y, \alpha)=M_p(y, \alpha)e^{-\alpha\mu(\alpha)y}.\ee
We obtain the relation
\be\label{evansfinal}\ba{ll}Ev(\alpha, \beta, \gamma)&=e^{-(2+\frac{\mu(\alpha)}{\beta})t}t^{\frac{1}{2\nu}}\exp(-\int_0^{-\frac{t_0}{\alpha\beta}}\alpha\gamma\xi(y')dy')\cr
&[z_1(R_1-\frac{\xi}{\beta^2}L_1-L_2)+\xi z_2 (R_2+\frac{\gamma}{\beta}L_1-L_4)\cr
&+\xi z_3(R_3+\frac{\gamma}{\beta}L_2-\frac{\xi}{\beta^2}L_4)+\xi z_4(R_4-\frac{\gamma}{\beta}L_3-\frac{\xi}{\beta^2}L_5-L_6)\cr
&+\xi\frac{m_1L_1+m_2L_2+m_3L_3}{1-\xi}\cr
&+\xi^2\frac{m_4L_4+m_5L_5+m_6L_6}{1-\xi}]\frac{1-\xi(0)}{\xi(0)}.\ea\ee
In  what follows we shall describe the functions $R$, $L$, and $z$, $m$.
\section{Calculus of the solution in the overdense region}
\label{over}
The aim of this section is to obtain $w_+^{(2)}$, which is the unique solution of the system (\ref{s2}) under the condition (\ref{conditionunicite2}) in the region $[y_0, +\infty[$ for all $y_0$. The first idea would be to try to apply the Gap Lemma. However, the eigenvalue of smallest real partof $-M_0^{(2)}(1)$ is $\lambda_{min}=\lambda_-(1)-\alpha\beta$, and the next eigenvalue is $\lambda^*=\lambda_-(1)+\alpha\gamma$. As $\lambda^*-\lambda_{min}=\alpha(\gamma + \beta)$, the difference between two eigenvalues of $M_0^{(2)}(1)$ is not uniformly bounded below for $\alpha\in ]0,\alpha_0[$. Hence the hypotheseses of the Gap lemma theorem are not fulfilled.\\ 
\paragraph{Remark} As the coefficients of the differential system behave as $1-e^{-y}$ when $y$ goes to $+\infty$, one may obtain $Y^{(2)}$ through a Volterra expansion as

$$Y^{(2)}(y)=Ae^{(\lambda_{ -}(1)-\alpha\beta) y}+ e^{-y}r(y, \alpha, \beta, \gamma)$$
$r$ being a remainder term. We prove that the result that we obtain is analytic in $\alpha$ for $\beta$, $\gamma$ in a certain compact set. Once we know that, we identify $w_+^{(2)}$ (deduced from $Y^{(2)}$ and obtained also through a Volterra expansion) and we express it with an expansion in powers of $\alpha$, which is bounded by a geometric series for $\xi(y)\in [\xi_0, 1]$.\\
We then plug this developement in the differential system and we identify the terms. We check that the coefficient of $\alpha^j$ in the expansion in $\alpha$ of the solution is of the form ${1\over \xi^{\nu j+2}}X_j(\xi)$, the radius of convergence of the series depends on $\xi(y_0)$. The aim of the explicit calculus is to deduce a new region of convergence of the series using the behavior in ${\alpha^j\over \xi^{\nu j+2}}X_j(\xi)$. In this new region of convergence, we have a converging series, which coincides with the original one for $\xi\in [\xi(y_0, 1[$. Hence it is the extension of our solution in the new region of convergence. This is an important feature, because it helps to have an overlapping region between the most decreasing solution at $-\infty$ and the most decreasing solution at $+\infty$.\\ We introduce $z_j$ and $m_p$ through
$$T^{(2)}w_+^{(2)}(y)=\alpha\beta(\sum_{j=1}^4z_jf_j+\sum_{p=1}^6m_pg_p)e^{\alpha\mu(\alpha)y}$$
where $\mu(\alpha)$ has been given in (\ref{defimu}). We recall that $z_j\rightarrow W_{+, j}$ and $m_j\rightarrow 0$ when $y\rightarrow +\infty$.\\
The following Theorem summarizes the results of this section. Let us introduce $\xi_0>0$.

\begin{theo}
\label{theozm}
\begin{enumerate}
\item The functions $z_j$ and $m_j$ have a normally convergent expansion in powers of $\alpha$ of the form:

$$\ba{l}z_1(y, \alpha)=\sum_{j\geq 1} {\alpha^j\over \xi^{\nu j}}a_{1, j}(\xi)+W_{+, 1}\cr
z_k(y, \alpha)=  \sum_{j\geq 1} {\alpha^j\over \xi^{\nu j+1}}a_{k, j}(\xi)+W_{+, k}, k=2, 3, 4\cr
m_l(y, \alpha)=\sum_{j\geq 1} {\alpha^j\over \xi^{\nu j+1}}b_{l, j}(\xi), l=1,2,3\cr
m_p(y, \alpha)=\sum_{j\geq 1} {\alpha^j\over \xi^{\nu j+2}}b_{p, j}(\xi), p=4,5,6.\ea$$
We introduce $\delta_1=0$, $\delta_2=\delta_3=\delta_4=d_1=d_2=d_3=1$, $d_4=d_5=d_6=2$.
Let $K$ be a compact subset of $\RR_+^*\times \RR$. There exists $R>0$ and $\alpha_1(\xi_0)$ such that for $\beta, \gamma$ in $K$ and $\alpha<\alpha_1(\xi_0)$ we have for $\xi\in [\xi_0, 1]$

$$\ba{l}|{d\over d\xi}({a_{k,j}(\xi)\over \xi^{\nu j+d_k}})|\leq \frac{R^j}{\xi^\nu j + d_k+1}\cr
|a_{k, j}(\xi)|\leq R^j(1-\xi),\cr
|{d\over d\xi}({b_{l,j}(\xi)\over \xi^{\nu j+\delta_l}})|\leq {R^j\over \xi^{\nu j + \delta_l+1}},\cr
|b_{l, j}(\xi)|\leq R^j(1-\xi).\ea$$
\end{enumerate}
\end{theo}
Assume this theorem is proven. We show that, for all $\xi_0$, there exists $\alpha_0(\xi_0)$ such that, for $\alpha<\alpha_0(\xi_0)$ the power series is convergent for $\xi\in [\xi_0, 1]$. Consider now ${\alpha\over \xi^{\nu}}R<1$ and ${\alpha\over \xi^{\nu}}M<1$. The power series $\sum_{j\geq 1} a_{p, j}(\xi){\alpha^j\over \xi^{\nu j + \delta_p}}$ defines an analytic function which is the analytic extension of the sum of the normally convergent series $\sum_{j\geq 1}{\alpha^j\over \xi^{\nu j + d_c}}a_{p, j}(\xi)$ for $\alpha<\alpha_1(\xi_0)$ and $\xi\in [\xi_0, 1]$. A similar result holds for the series $\sum_{j\geq 1}b_{l, j}(\xi)\frac{\alpha^j}{\xi^{\nu j + \delta_l}}$.

The proof of item 1) of the theorem is in Annex \ref{annexe1}. We calculate in the first subsection of the present section the first order terms of $z$ and $m$. \\
In the second subsection, we prove by recurrence the structure of the $j-$th term of the expansion of $z$ and $m$, which is a consequence of the structure of the system. We deduce the precise estimates wich help us to extend the expansion.

\subsection{A new formulation of the system}
We introduce $Z_j, 1\leq j\leq 4$ and $m_p$, $1\leq p\leq 6$ given by

\be\label{formulew}{1\over \alpha \beta}T^{(2)}w_+^{(2)}(y, \alpha)=y^{(2)}_+(y, \alpha)=\sum_{j=1}^4 Z_jf_j+\sum_{p=1}^6M_pg_p=e^{\alpha\mu(\alpha)y}(\sum_{j=1}^4 z_jf_j+\sum_{p=1}^6m_pg_p).\ee
The relation
$w_+^{(2)}e^{y}e^{-\alpha\mu(\alpha)y}\rightarrow F_-(1)\wedge E_{a, -}(1)=W_+$ imply that $z_j\rightarrow W_{+, j}$ because ${\xi\over (1-\xi)e^y}\rightarrow 1$ when $y\rightarrow +\infty$.
The system on $Z, M$ writes
\be\bs\ba{l}{dZ\over dy} + {\alpha\over
  \xi}J_0Z+\alpha B(\xi)Z+{\alpha\beta\over 1-\xi}(L_0+\xi L_1)M=0\cr
{dM\over dy} +{\alpha\over \xi}K_0 M +\alpha D(\xi) M
  +\alpha{1-\xi\over \xi} C(\xi)Z=0.\ea\es\ee
 that is
 \be\label{S1}\bs\ba{ll}{1\over \alpha}{dZ_1\over dy}+&\beta \xi Z_3 -\gamma\xi Z_4 -{\beta\xi M_3\over 1-\xi}=0\cr
  {1\over \alpha}{dZ_2\over dy} + &(\gamma -{1\over \beta})Z_1+\gamma Z_2 -\beta Z_3 -{\xi\over \beta}Z_4 -{\beta (M_1+\xi M_5)\over 1-\xi}=0\cr
  {1\over \alpha}{dZ_3\over dy} + &{\beta\over \xi}Z_1 + \beta Z_2+\gamma(1-\xi)Z_3-\beta Z_4 -{\beta (M_2+\xi M_6)\over 1-\xi}=0\cr
  {1\over \alpha}{dZ_4\over dy} &- \beta Z_3+\gamma Z_4 +{\beta M_3\over 1-\xi}=0\cr
  {1\over \alpha}{dM_1\over dy} &-\gamma M_1 - \beta M_2 +{\xi\over \beta}M_3 + \xi (-\beta M_4 -\gamma M_5)\cr
  & + (1-\xi)[({\gamma^2\over \beta} - {\gamma\over \beta^2} -{\beta \over \xi})Z_1 + {\gamma^2-\beta^2\over \beta}Z_2]=0\cr
  {1\over \alpha}{dM_2\over dy} &+ \beta M_1 -\gamma(1+\xi)M_2 + \beta M_3 -\gamma\xi M_6\cr
  & + (1-\xi)[({1\over \beta}+\gamma{1-\xi\over \xi})Z_1 +{\gamma^2-\beta^2\over \beta}Z_3]=0\cr
  {1\over \alpha}{dM_3\over dy} &+ \beta M_2 -\gamma M_3 + \beta\xi M_6 + (1-\xi){\beta^2-\gamma^2\over \beta}Z_4=0\cr
  {1\over \alpha}{dM_4\over dy}&-{\beta\over \xi}M_1 + (\gamma -{1\over \beta})M_2 -\gamma\xi M_4 + \beta M_5 -{\xi\over \beta}M_6\cr
  &+(1-\xi)[Z_2({1\over \beta}+\gamma{1-\xi\over \xi}) + Z_3({\beta\over \xi}+{\gamma\over \beta^2}-{\gamma^2\over \beta})]=0\cr
  {1\over \alpha}{dM_5\over dy}&+(\gamma -{1\over \beta})M_3 + \beta M_4 -\beta M_6 +(1-\xi)Z_4(-{\beta\over \xi}+{\gamma\over \beta}(\gamma -{1\over \beta}))=0\cr
  {1\over \alpha}{dM_6\over dy} &+ {\beta\over \xi}M_3 + \beta M_5 -\gamma\xi M_6 + (1-\xi)Z_4({1\over \beta}+\gamma{1-\xi\over \xi})=0.\ea\es\ee
This system writes

\be\label{ou1}{1\over \alpha}{dU\over dy}+\beta K(r, \xi(y), \beta)U=0\ee
where $U^t=(Z_1, Z_2, Z_3, Z_4, m_1, m_2, m_3, m_4, m_5, m_6)$.
The system on $(z, m)$ is deduced from (\ref{ou1}) by replacing the matrix $\beta K(r, \xi(y))$ by $\beta K(r, \xi(y))+\mu(\alpha)I$.
We verify that $m_p\rightarrow 0$ when $y\rightarrow +\infty$. In what follows, we describe the analytic expansion of the solution in $\alpha$ when $\alpha$ is in a neighborhood of 0.

\subsection{First terms of the expansion in $\alpha$}
From
\be\label{limites} z_j\rightarrow W_{+, j}, m_j\rightarrow 0\ee
we obtain $z_j^0=W_{+, j}, m_p^0=0$, that is $z_1^0=\beta, z_2^0=\gamma - \beta, z_3^0=\beta, z_4^0=-\beta$.
We replace in the system these relations to obtain the system on the term $z_j^1$ and $m_j^1$. We get

$$\bs\ba{l}{dz_1^1\over dy}+\mu(0)z_1^0+\beta\xi z_3^0-\gamma\xi z_4^0=0\cr
{dz_2^1\over dy} +\mu(0)z_2^0+ (\gamma -{1\over \beta})z^0_1+\gamma z^0_2 -\beta z^0_3 -{\xi\over \beta}z^0_4 =0\cr
  {dz^1_3\over dy}  +\mu(0)z_3^0+ {\beta\over \xi}z^0_1 + \beta z^0_2+\gamma(1-\xi)z^0_3-\beta z^0_4 =0\cr
  {dz^1_4\over dy}  +\mu(0)z_4^0- \beta z^0_3+\gamma z^0_4 =0\cr
  {dm^1_1\over dy} +(1-\xi)[({\gamma^2\over \beta} - {\gamma\over \beta^2} -{\beta \over \xi})z^0_1 + {\gamma^2-\beta^2\over \beta}z^0_2]=0\cr
  {dm^1_2\over dy}  + (1-\xi)[({1\over \beta}+\gamma{1-\xi\over \xi})z^0_1 +{\gamma^2-\beta^2\over \beta}z^0_3]=0\cr
  {dm^1_3\over dy}  + (1-\xi){\beta^2-\gamma^2\over \beta}z^0_4=0\cr
  {dm^1_4\over dy}+(1-\xi)[z^0_2({1\over \beta}+\gamma{1-\xi\over \xi}) + z^0_3({\beta\over \xi}+{\gamma\over \beta^2}-{\gamma^2\over \beta})]=0\cr
  {dm^1_5\over dy} +(1-\xi)z^0_4(-{\beta\over \xi}+{\gamma\over \beta}(\gamma -{1\over \beta}))=0\cr
  {dm^1_6\over dy}  + (1-\xi)z^0_4({1\over \beta}+\gamma{1-\xi\over \xi})=0\ea\es$$
As $\mu(0)=-\beta -\gamma$ we deduce
$$\bs\ba{l}{dz_1^1\over dy}-\beta(\beta+\gamma)(1-\xi)=0\cr
{dz_2^1\over dy} +\xi-1 =0\cr
  {dz^1_3\over dy}  +{\beta^2(1-\xi)\over \xi} +\gamma\beta (1-\xi) =0\cr
  {dz^1_4\over dy}  =0\cr
  {dm^1_1\over dy} +(1-\xi)[({\gamma^2\over \beta} - {\gamma\over \beta^2} -{\beta \over \xi})\beta + {\gamma^2-\beta^2\over \beta}(\gamma - \beta)]=0\cr
  {dm^1_2\over dy}  + (1-\xi)[({1\over \beta}+\gamma{1-\xi\over \xi})\beta +{\gamma^2-\beta^2\over \beta}\beta]=0\cr
  {dm^1_3\over dy}  - (1-\xi){\beta^2-\gamma^2\over \beta}\beta=0\cr
  {dm^1_4\over dy}+(1-\xi)[(\gamma - \beta)({1\over \beta}+\gamma{1-\xi\over \xi}) + \beta({\beta\over \xi}+{\gamma\over \beta^2}-{\gamma^2\over \beta})]=0\cr
  {dm^1_5\over dy} -\beta(1-\xi)(-{\beta\over \xi}+{\gamma\over \beta}(\gamma -{1\over \beta}))=0\cr
  {dm^1_6\over dy}  -\beta (1-\xi)({1\over \beta}+\gamma{1-\xi\over \xi})=0\ea\es$$
We obtain
$$\bs\ba{l}z_1^1= \beta(\beta+\gamma){\xi^{\nu}-1\over \nu \xi^{\nu}}\cr
z_2^1= {\xi^{\nu}-1\over \nu \xi^{\nu}}\cr
z_3^1=\gamma\beta{1-\xi^{\nu}\over \nu \xi^{\nu}}+\frac{\beta^2}{\nu+1}{1-\xi^{\nu+1}\over \xi^{\nu+1}}\cr
z_4^1=0\ea\es$$
and the following system on $m_j^1$:

$$\bs\ba{l}
  {dm^1_1\over dy} +{{\dot \xi}\over \xi^{\nu+2}}[(\gamma^2\xi - {\gamma\over \beta}\xi -\beta^2) + {\gamma^2-\beta^2\over \beta}(\gamma - \beta)\xi]=0\cr
  {dm^1_2\over dy}  + {{\dot \xi}\over \xi^{\nu+2}}[\xi +\gamma\beta (1-\xi) +(\gamma^2-\beta^2)\xi]=0\cr
  {dm^1_3\over dy}  - {{\dot \xi}\over \xi^{\nu+1}}(\beta^2-\gamma^2)=0\cr
  {dm^1_4\over dy}+{{\dot \xi}\over \xi^{\nu+2}}[(\gamma - \beta)({\xi\over \beta}+\gamma(1-\xi)) + (\beta^2+{\gamma\over \beta}\xi-\gamma^2\xi)]=0\cr
  {dm^1_5\over dy} -{{\dot \xi}\over \xi^{\nu+2}}(-\beta^2+\gamma\xi(\gamma -{1\over \beta}))=0\cr
  {dm^1_6\over dy}  -{{\dot \xi}\over \xi^{\nu+2}}(\xi+\gamma\beta (1-\xi))=0\ea\es$$
Hence we obtain the expansion of the solution for $\xi\geq \xi_0$

$$\bs\ba{l}z_1=\beta -{\alpha\over \nu \xi^{\nu}}(1-\xi^{\nu})\beta(\beta+\gamma)+O(\alpha^2)\cr
z_2=\gamma -\beta - {\alpha\over \nu \xi^{\nu}}(1-\xi^{\nu})+O(\alpha^2)\cr
z_3=\beta +{\alpha\over (\nu+1)\xi^{\nu+1}}\beta^2(1-\xi^{\nu+1})+{\alpha\over \nu \xi^{\nu}}\gamma\beta(1-\xi^{\nu})+O(\alpha^2)\cr
z_4=-\beta+O(\alpha^2)\ea\es$$
The next subsection is dedicated to the precise study of the behavior of $w_+^{(2)}$, which depends on inverse powers of $\xi$ and cannot be extended directly to $\xi\rightarrow 0$.

\subsection{Uniform estimates of the solution $w_+^{(2)}$ for $\xi\in [{\zeta_0\over \alpha^{1\over \nu}}, 1]$}
\label{suniform}
We write
\be\label{calculz}z=W_{+, 1}f_1+W_{+, 2}fÃ_2+W_{+, 3}f_3+W_{+, 4}f_4+u.\ee
We introduce the constant matrices $J_0, L_0$ and $K_0$ such that $J=J_0+\xi B(\xi), L=L_0+\xi L_1, K=K_0+\xi D(\xi)$. The system (\ref{sl}) yields

\be\label{systemeavecmu}\bs\ba{l}{dz\over dy}+\alpha\mu(\alpha)z + {\alpha\over \xi}Jz+{\alpha\beta\over 1-\xi}Lm=0\cr
{dm\over dy}+\alpha\mu(\alpha)m + {\alpha\over \xi}Km + \alpha{1-\xi\over \xi}Cz=0\ea\es\ee

The aim of this paragraph is to find a simpler formulation for the unique solution of this system.

As the solution going to 0 at infinity of ${df\over dy}={1-\xi\over \xi}$ is $f(y)={1\over \nu+1}(1-{1\over \xi^{\nu+1}})$, direct estimates of the behavior in $\xi$ of a coefficient $u_j$ or $m_j$ of (\ref{systemeavecmu}) lead to a multiplying factor of the form

$${\alpha\over \xi^{\nu+1}}.$$
The behavior of $\alpha^{j+1}u_{j+1}$ or of $\alpha^{j+1}m_{j+1}$ in $\xi$ is then given by ${\alpha^{j+1}\over \xi^{(\nu+1)j}}$ for the next coefficient.
However, the structure of the system allows us to obtain a lower inverse power of $\xi^\nu$ in the expansion in $\alpha$. For this purpose, we introduce the new unknowns $a_{p, j}$, $b_{q, j}$ such that

\be\label{zm}\bs\ba{l}z_1=W_{+, 1}+\sum_{j=1}^N{\alpha^j\over \xi^{\nu j}}a_{1, j}(\xi)+z_1^{N+1}(\alpha, \xi)\cr
z_p=W_{+, p}+\sum_{j= 1}^N{\alpha^j\over \xi^{\nu j+1}}a_{p, j}(\xi)+z_p^{N+1}(\alpha, \xi), p=2,3,4\cr
m_k=\sum_{j=1}^N{\alpha^j\over \xi^{\nu j+1}}b_{k, j}(\xi)+m_k^{N+1}(\alpha, \xi), k=1,2,3\cr
m_l=\sum_{j= 1}^N{\alpha^j\over \xi^{\nu j+2}}b_{l, j}(\xi)+m_l^{N+1}(\alpha, \xi), l=4,5,6\ea\es\ee
The previous quick analysis would suggest that $a_{p, j+1}$ is of order ${1\over \xi}$ when every $a_{k, j}, b_{l, j}$ is bounded when $\xi$ goes to 0. This is not the case, and the crucial equality states as follows. We introduce the diagonal matrix $T$ such that

$$T\left(\ba{c}a_1\cr a_2\cr a_3\cr a_4\cr b_1\cr b_2\cr b_3\cr b_4\cr b_5\cr b_6\ea\right)= \left(\ba{c}a_1\cr \xi^{-1}a_2\cr \xi^{-1}a_3\cr \xi^{-1}a_4\cr \xi^{-1}b_1\cr \xi^{-1}b_2\cr \xi^{-1}b_3\cr \xi^{-2}b_4\cr \xi^{-2}b_5\cr \xi^{-2}b_6\ea\right).$$
There exists a matrix ${\hat C}(\xi, \beta, \gamma)$, polynomial in $\xi$, such that

\be\left(\ba{cc}{1\over \xi}J&{\beta\over 1-\xi}L\cr
{1-\xi\over \xi}C&{1\over \xi}K\ea\right)T= T{\hat C}.\ee
We shall make use of the following fundamental Lemma, noting that at each step we solve an equation of the form
$${df\over dy}={A(\xi)(1-\xi)\over \xi^\alpha}$$
where $\alpha=\nu j + d$, $d=0, 1, 2$.
\begin{lemma}
\label{primitive}
The unique solution going to 0 when $\xi$ goes to 1 of

$${df\over dy} = {A(\xi)(1-\xi)\over \xi^\alpha}$$
is $f(y)=\int_\xi^1{A(\eta)\over \eta^{\alpha+\nu+1}}d\eta$. We have the estimate

$$|f(y)|\leq {1-\xi\over \xi^{\nu+\alpha}}||A||_\infty .$$
\end{lemma}
\paragraph{Proof} From ${1\over \xi^{\nu + 1 + \alpha}}={d\over d\xi}({1\over \nu + \alpha}
(1-{1\over \xi^{\nu+\alpha}}))$, we deduce
$$|f(y)|\leq ||A||_{\infty}\int_{\xi}^1{d\over d\xi}({1\over \nu + \alpha}
(1-{1\over \eta^{\nu+\alpha}}))d\eta=
||A||_{\infty}{1-\xi^{\nu + \alpha}\over (\nu+\alpha)\xi^{\nu + \alpha}}.$$
The equality ${\xi^{\nu + \alpha} - 1\over \xi - 1}=\int_0^1(\nu + \alpha)(1+t(\xi-1))^{\nu+\alpha -1}dt$ implies
\be\label{I}|\xi^{\nu + \alpha} - 1|\leq (\nu + \alpha)(1-\xi).\ee
hence the lemma.
The indices defined in Theorem \ref{theozm} will express the weight of each coordinate of the vector $U$ defined in (\ref{ou1})
$$T\left[{1\over \xi^{\nu j}}\left(\ba{c}a_j\cr b_j\ea\right)\right]=\left(\ba{c}{a_{p, j}\over \xi^{\nu j + \delta_j}}\cr {b_{q, j}\over \xi^{\nu j + d_q}}\ea\right).$$
In the system (\ref{zm}), write

$$\left(\ba{c}z\cr m\ea\right)=\left(\ba{c}W_+\cr 0\ea\right)+\sum_{j\geq 1}{\alpha^j\over \xi^{\nu j}}T\left(\ba{c}a_j\cr b_j\ea\right)$$
we get, for the term in $\alpha^{j+1}$:

$${d\over dy}({1\over \xi^{\nu(j+1)}}T\left(\ba{c}a_{j+1}\cr b_{j+1}\ea\right))+\mu_j\left(\ba{c}W_+\cr 0\ea\right) +\sum_{l=1}^j{\mu_{j-l}\over \xi^{\nu l}}T\left(\ba{c}a_l\cr b_l\ea\right)+{1\over \xi^{\nu j}}TC\left(\ba{c}a_j\cr b_j\ea\right)=0.$$
This system of equations becomes, for $j\geq 2$:
\be\bs\ba{l}{d\over dy}({a_{p, j+1}\over \xi^{\delta_p+\nu(j+1)}})+\frac{(\sum_{l=0}^j\mu_la_{p, j-l}\xi^{\nu l}) + (C_{11}a_j+C_{12}{b_j\over 1-\xi}))_p}{\xi^{\delta_p+\nu j}}+\mu_{j-1}W_{+, p}=0\cr
{d\over dy}({b_{p, j+1}\over \xi^{d_p+\nu(j+1)}}) + \frac{(\sum_{l=0}^j\mu_lb_{p, j-l}) + C_{21}(1-\xi)a_j+C_{22}b_j)}{\xi^{\delta_p+\nu (j+1)}}=0\ea\es\label{systemerecurrent}\ee
and the equality for $j=1$
\be\bs\ba{l}{d\over dy}({a_{p, j+1}\over \xi^{\delta_p+\nu(j+1)}})+\frac{(\sum_{l=0}^j\mu_la_{p, j-l}\xi^{\nu l}) + (C_{11}a_j+C_{12}{b_j\over 1-\xi}))_p}{\xi^{\delta_p+\nu j}}+ (1-\xi)h_p=0\cr
{d\over dy}({b_{p, j+1}\over \xi^{d_p+\nu(j+1)}}) + \frac{ (\sum_{l=0}^j\mu_lb_{p, j-l}) + C_{21}(1-\xi)a_j+C_{22}b_j)}{\xi^{d_p+\nu(j+1)}}+(1-\xi)h_{4+p}=0.\ea\es\ee
We have the identity $C_{12}(1)b'_j(1)=h_j(1)+\mu_{j-1}W_{+, j}$, which is necessary to obtain that the source term in the equation on $a_{p, j+1}$ vanishes at $\xi=1$.
\subsection{Behavior of the terms of the expansion}
\label{sprecisees}
The regularity of the quantities $a_{p,j}$, $b_{q,j}$ is given by the following proposition, which gives precise estimates on the functions provided that $\beta$ and $\gamma$ stay in a compact set:
\begin{prop}
\label{inegalitesprecisees}
Assume that $\beta, \gamma$ are in a compact set $K$ of $\RR_+^*\times \CC$, namely $$\beta_0\leq \beta\leq {1\over \beta_0}, |\gamma|\leq {1\over \beta_0},\Re\gamma\geq 0.$$
Assume that $\xi_0\in ]0, 1[$ is given and that there exists $\alpha_0(\xi_0)$ such that the analytic expansion of the solution of (\ref{systemeavecmu}) is valid for $\alpha<\alpha_0(\xi_0)$ and $\xi(y)\in [\xi_0, 1[$.

 There exists $R$ and $M$ depending only of $\alpha_0$ and $\beta_0$ such that, forall $p=1, 2, 3, 4$, forall $q=1, ..., 6$, for all $j\geq 1$ we have the estimates
$$|a_{p,j}(\xi)|\leq R^j(1-\xi)$$
$$|b_{q,j}(\xi)|\leq R^j(1-\xi)$$
If we introduce $\alpha_j= a_j(0)$ and $\beta_j = b_j(0)$ we have
$$|a_j(\xi)-\alpha_j(1-\xi)|\leq M^j\xi(1-\xi)$$
$$|b_j(\xi)-\beta_j(1-\xi)|\leq M^j\xi(1-\xi).$$
\end{prop}
From this proposition we deduce the

\begin{prop}
Under the same hypothesis as the proposition \ref{inegalitesprecisees},

i) there exists $\alpha_1(\xi_0)$ such that, for $\alpha< \alpha_1(\xi_0)$, the functions
$$\ba{l}A_p(\zeta, \alpha)=\sum_{j=1}^{\infty}a_{p,j}(({\alpha\over \zeta})^{1\over \nu})\zeta^j\cr
B_q(\zeta, \alpha)=\sum_{j=1}^{\infty}b_{q,j}(({\alpha\over \zeta})^{1\over \nu})\zeta^j\ea$$
are analytic for $\zeta<{1\over R}$.

 2) For $y$ such that

$$\alpha^{1\over \nu}R^{1\over \nu}<\xi(y)\leq 1$$
the functions

$$u_p(y, \alpha) = W_{+, p} + \xi^{-\delta_p}A_p({\alpha\over (\xi(y))^\nu}, \alpha), v_q(y, \alpha)=\xi^{-d_q}B_q({\alpha\over (\xi(y))^\nu}, \alpha)$$ are solution of the system (\ref{systemeavecmu}) and extend the gap lemma solution.

3) Introduce

$${\tilde A}_p(\zeta, \alpha)=\sum_{j=1}^{\infty}\zeta^ja_{p,j}(({\alpha\over \zeta})^{1\over \nu})-\alpha_{p,j}(1-({\alpha\over \zeta})^{1\over \nu})\zeta^j.$$
$${\tilde B}_q(\zeta, \alpha)=\sum_{j=1}^{\infty}\zeta^jb_{q,j}(({\alpha\over \zeta})^{1\over \nu})-\beta_{q,j}(1-({\alpha\over \zeta})^{1\over \nu})\zeta^j.$$
The functions ${\tilde A}_p$ and ${\tilde B}_q$ are analytic for $\zeta<{1\over M}$, and we have the inequalities for $\zeta<{1\over M}$
$$|{\tilde A}_p(\zeta, \alpha)|\leq ({\alpha\over \zeta})^{1\over \nu}{\zeta\over 1-M\zeta}.$$
$$|{\tilde B}_q(\zeta, \alpha)|\leq ({\alpha\over \zeta})^{1\over \nu}{\zeta\over 1-M\zeta}=\alpha^{1\over \nu}{\zeta^{1-{1\over \nu}}\over 1-M\zeta}.$$
\label{serieprolongee}
\end{prop}
We deduce, for $\alpha\leq\alpha_1(\xi_0)$ and $\zeta\leq {1\over M}$, the equalities

$$\bs\ba{l}U_1(\zeta, \alpha)=W_{+, 1}+\sum_{j=1}^{\infty}\alpha_{1, j}\zeta^j+\alpha^{1\over \nu}R_1(\zeta, \alpha)\cr
\alpha^{1\over \nu}U_p(\zeta, \alpha)=\zeta^{1\over \nu}\sum_{j=1}^{\infty}\alpha_{p, j}\zeta^j+\alpha^{1\over \nu}R_p(\zeta, \alpha), p=2,3,4\cr
\alpha^{1\over \nu}V_q(\zeta, \alpha)=\zeta^{1\over \nu}\sum_{j=1}^{\infty}\beta_{q, j}\zeta^j+\alpha^{1\over \nu}S_q(\zeta, \alpha), q=1,2,3\cr
\alpha^{2\over \nu}V_q(\zeta, \alpha)=\zeta^{2\over \nu}\sum_{j=1}^{\infty}\beta_{q, j}\zeta^j+\alpha^{1\over \nu}S_q(\zeta, \alpha), q=4,5,6\ea\es$$
\paragraph{Proof of Proposition \ref{serieprolongee}}

Let $\xi_0\in ]0, 1[$ be given. There exists $\alpha_0(\xi_0)$ such that, for $\alpha\leq\alpha_0(\xi_0)$, the solution satisfies the conclusions of the gap lemma (which means that the solution of the system (\ref{systemeavecmu}) is analytic in the region $\alpha\leq \alpha_0(\xi_0)$ for $\xi\in [\xi_0, 1]$).

Moreover, from proposition \ref{inegalitesprecisees}, for all $\xi_0$ there exists $R$ depending only on $\beta_0$ such that for $\zeta<{1\over R}$ the functions $A_p(\zeta, \alpha)$ and the functions $B_q(\zeta, \alpha)$ given in Proposition \ref{serieprolongee}  are analytic through their expansion in $\alpha$ for $\zeta<{1\over R}$.\\
Introduce $\alpha_1(\xi_0)=\mbox{min}(\alpha_0(\xi_0), {\xi_0^\nu\over R})$. For $\alpha<\alpha_1(\xi_0)$ and $\zeta<{1\over R}$, $\xi=({\alpha\over \zeta})^{1\over \nu}\geq \xi_0$. This means that

$$(w_{+, p}+({\zeta\over \alpha})^{\delta_p\over \nu}A_p(({\alpha\over \zeta})^{1\over \nu}, \alpha), ({\zeta\over \alpha})^{d_q\over \nu}B_q(({\alpha\over \zeta})^{1\over \nu}, \alpha))$$
is solution of the system (\ref{systemeavecmu}) for $\alpha\leq\alpha_0(\xi_0)$ when $y$ is given by $\xi(y)=({\alpha\over \zeta})^{1\over \nu}$ by construction of the analytic solution given by the gap lemma.

By uniqueness of the solution which is analytic in $\alpha$, we check that this function is also solution of the system (\ref{systemeavecmu}) for $\alpha<\alpha_1(\xi_0)$ and $\xi(y)\in [(\alpha R)^{1\over \nu}, \xi_0]$, because the analytic expansion defining the solution in the set-up of the gap lemma can be rearranged and the remainder term is regular enough (and uniformly bounded), and because the two solutions are equal at the point ${\tilde y}$ such that $\xi({\tilde y})={1+\xi_0\over 2}$.

Hence we extended the solution for $\xi$ in the interval $[(\alpha R)^{1\over \nu}, 1]$. Proposition \ref{serieprolongee} is thus a consequence of Proposition \ref{inegalitesprecisees}.

\paragraph{Proof of Proposition \ref{inegalitesprecisees}}
We prove by recurrence the inequalities (\ref{a1}), (\ref{a2}), (\ref{a3}) below:

\be\label{a1} |{d\over d\xi}({a_{p, j}\over \xi^{\delta_p+\nu j}})|\leq
{R^j\over \xi^{\nu j + 1 + \delta_p}}.\ee

\be\label{a2} |{d\over d\xi}({b_{q, j}\over \xi^{d_q+\nu j}})|\leq
{R^j\over \xi^{\nu j + 1 + d_q}}.\ee

\be\label{a3} |{d\over d\xi}({b_{q, j}\over\xi^{d_q+\nu j}})(\xi)
-{d\over d\xi}({b_{q, j}\over\xi^{d_q+\nu j}})(1)|\leq {R^j(1-\xi)\over \xi^{\nu j + 1 + d_q}}.\ee
We assume that $(\beta, \gamma)\in K\subset[\beta_0, {1\over \beta_0}]\times \{\gamma, |\gamma|\leq \beta_0^{-1}, \Re\gamma\geq 0\}$.
\paragraph{First step:}
From the identity
$${a_{p, j}(\xi)\over \xi-1}= \xi^{\delta_p+\nu j}\int_1^\xi
{d\over d\xi}({a_{p, j}\over \xi^{\delta_p+\nu j}})
(1+s(\xi-1))ds$$
and from the same identity on $b_j$, relying on the fact that $a_j$ and $b_j$ are 0 at $\xi=1$ we deduce the inequality

$$|a_{p, j}(\xi)|\leq R^j\xi^{\delta_p+ \nu j}
\int_\xi^1{(1-\xi)ds\over (1+s(\xi-1))^{\nu j + 1 + \delta_p}}.$$
We use
 (\ref{I}) to obtain

$$|a_{p, j}(\xi)|\leq R^j(1-\xi),$$
or
\be\label{a4}|{a_{p, j}(\xi)\over 1-\xi}|\leq R^j\ee
\be\label{a5}|{b_{q, j}(\xi)\over 1-\xi}|\leq R^j\ee
and of course
\be\label{a6}|{db_{p, j}\over d\xi}(1)|=|{d\over d\xi}({b_{p, j}\over \xi^{d_p+ \nu j}})(1)|\leq R^j.\ee

We assume that (\ref{a1}), (\ref{a2}), (\ref{a3}) are true for $l\leq j$. We deduce the inequalities (\ref{a4}), (\ref{a5}), (\ref{a6}).

The equation on $a_{p, j+1}$ rewrites

$$\ba{l}\xi^{\nu+1}(1-\xi){d\over d\xi}({a_{p, j+1}\over \xi^{\delta_p+\nu(j+1)}}) + {1\over \xi^{\nu j + \delta_p}}[C_{11}(\xi)a_j(\xi)+\sum_{l=0}^j\mu_l\xi^{\nu l} a_{p, j-l}\cr
+(C_{12}{b_j\over 1-\xi}+C_{12}(1)b'_j(1))].\ea$$

Recall that we have the equality

$$-\mu(\alpha)= \beta+2{\gamma+\alpha \beta^2\over 1+\sqrt{1+4\alpha(\gamma+\alpha\beta^2)}}$$
hence for $(\beta, \gamma)$ in the compact $K$, the function $\mu(\alpha)$ admits a DSE at $\alpha=0$, of radius of convergence greater than $\theta_0= \mbox{min}(1, {\beta_0\over 8})$. Moreover, denoting by $C_0=\sum_{l=0}^{\infty}|\mu_l|\theta_0^l$, we have

$$\sum_{l=0}^{\infty}|\mu_l|\theta^l\leq C_0.$$
Using the norm of the matrices $C_{11}+\mu_0I$, $C_{12}$, $C_{21}$ and $C_{22}+\mu_0I$ (with $0\leq \xi\leq 1$) and the inequalities (\ref{a4}), (\ref{a5}), (\ref{a6}) we obtain

$$\ba{ll}|{d\over d\xi}({a_{p, j+1}\over \xi^{\delta_p+\nu(j+1)}})|\leq
&{1\over \xi^{\nu(j+1) + 1+\delta_p}}[\sum_{l=1}^j \mu_lR^{j-l}\xi^{\nu l}
+ (3|\gamma| + 6\beta + {2|\gamma|^2+3\over \beta}+{|\gamma|\over \beta^2})R^j]\cr
&\leq {R^j\over \xi^{\nu(j+1) + 1+\delta_p}}[C_0+3|\gamma| + 6\beta + {2|\gamma|^2+3\over \beta}+{|\gamma|\over \beta^2}] \ea$$
as soon as \be\label{inegR}R^{-1}\leq \mbox{min}(1, {\beta_0\over 8}).\ee

We have the same type of estimates for $b_{q, j+1}$:

$$|{d\over d\xi}({b_{q, j+1}\over \xi^{d_q+\nu(j+1)}})| \leq {R^j\over \xi^{\nu(j+1)+1+d_q}}[C_0+3|\gamma| + 6\beta + {2|\gamma|^2+3\over \beta}+{|\gamma|\over \beta^2}].$$
Hence there exists a constant $D(\beta_0)$ depending only on the compact set $K$, such that for $R$ satisfying (\ref{inegR}) and assuming (\ref{a1}), (\ref{a2}), (\ref{a3}) at the order $j$, we obtain the estimates

$$|{d\over d\xi}({a_{p, j+1}\over \xi^{\delta_p+\nu (j+1)}})|\leq
{R^jD(\beta_0)\over \xi^{\nu (j+1) + 1 + \delta_p}}$$

$$ |{d\over d\xi}({b_{q, j+1}\over \xi^{d_q+\nu (j+1)}})|\leq
{R^jD(\beta_0)\over \xi^{\nu (j+1) + 1 + d_q}}.$$
The last estimate that we need is based on the difference of derivatives and we use the identity

$${f(\xi)\over \xi^{\alpha}(1-\xi)}+f'(1)={1\over \xi^{\alpha}}\int_0^1(f'(1)-f'(1+s(\xi-1)))ds -{f'(1)\over \xi^{\alpha}}(1-\xi^{\alpha})$$
for $f(1)=0$ and $f\in C^1$. Assume the inequalities $|f'(\xi)-f'(1)|\leq C{(1-\xi)\over \xi^{\beta+1}}$ and $|f'(1)|\leq C$. We obtain

$$|{f(\xi)\over \xi^{\alpha}(1-\xi)}+f'(1)|\leq {(1-\xi)C(1+\alpha)\over \xi^{\alpha+\beta}}.$$
We apply this estimate for $f(\xi)={b_{p,j- l}\over \xi^{d_p+\nu l}}$ and $\alpha=\nu+1$ ($\beta=d_p+\nu(j-l)$) to obtain the inequality
$$|{b_{p, j-l}\over (1-\xi)\xi^{\nu(j-l)+d_p+\nu + 1}}+b'_{p, j-l}(1)|\leq {1\over \xi^{\nu+1}}{R^{j-l}(1-\xi^{\nu(j-l)+d_p})\over (\nu(j-l)+d_p)\xi^{\nu(j-l)+d_p}}+R^{j-l}(\nu+1){1-\xi\over \xi^{\nu+1}}.$$
We thus deduce the estimate

$$|{d\over d\xi}({b_{q, j+1}\over\xi^{d_q+\nu (j+1)}})(\xi)
-{d\over d\xi}({b_{q, j+1}\over\xi^{d_q+\nu (j+1)}})(1)|\leq {R^jD(\beta_0)(1-\xi)\over \xi^{\nu (j+1) + 1 + d_q}}.$$

We use the estimate for $j=1$, for which there exists a constant $D_1(\beta_0)$ such that
$$|{d\over d\xi}({a_{p, 1}\over \xi^{\delta_p+\nu }})|\leq
{D_1(\beta_0)\over \xi^{\nu + 1 + \delta_p}}, |{d\over d\xi}({b_{q, 1}\over \xi^{d_q+\nu }})|\leq
{D_1(\beta_0)\over \xi^{\nu + 1 + d_q}}$$
$$|{d\over d\xi}({b_{q, 1}\over\xi^{d_q+\nu }})(\xi)
-{d\over d\xi}({b_{q, 1}\over\xi^{d_q+\nu }})(1)|\leq {D_1(\beta_0)(1-\xi)\over \xi^{\nu + 1 + d_q}}.$$
It is then enough to consider

$$R=\mbox{max}(1, {8\over \beta_0}, D(\beta_0), D_1(\beta_0))$$
to obtain the inequalities (\ref{a1}), (\ref{a2}), (\ref{a3}) for all $j$.
We write each term $a_{p, j}=\alpha_{p, j}+\xi c_{p, j}$, $b_{q, j}=\beta_{q, j}+\xi c_{q, j}$, and we have similar inequalities for the terms $c_{p, j}$ and $d_{q, j}$, with a coefficient $M$ depending only on the compact set $K$.
To obtain the estimate of the rest (versus the leading order term), we denote by $c_j$ and $d_j$ the functions such that $a_{j, p}(\xi)=(1-\xi)(\alpha_{j,p}+\xi c_{j, p}(\xi))$, $b_{q, j}(\xi)=(1-\xi)(\beta_{j, q}+\xi d_{q, j}(\xi))$.
We prove in a similar fashion the inequalities

\be\label{b1}|{d\over d\xi}({c_{p, j}(\xi)\over \xi^{\delta_p+\nu j}})|\leq {M^j\over \xi^{\nu j + \delta_p}}\ee
\be\label{b2}|{d\over d\xi}({d_{j, q}(\xi)\over \xi^{d_q}})|\leq {M^j\over \xi^{\nu j + d_q}}\ee

\be\label{b3}|{d\over d\xi}({d_{j, q}(\xi)\over \xi^{d_q}})(\xi) +d'_{j, q}(1)|\leq {M^j(1-\xi)\over \xi^{\nu j + d_q}}.\ee
This ends the proof of Proposition \ref{inegalitesprecisees}.

\section{The simplest discontinuity model}
\label{petitmodele}
Before studying the coupling between the hypergeometric region and the overdense region, we shall in this section study a simple model where the profile of density is
$$\xi(y)=\bs\ba{l}\xi(y_0), y\in ]-\infty, y_0]\cr
\xi(y), y\in ]y_0, +\infty[.\ea\es$$
It is a slightly better model than the discontinuity model (see \cite{Piriz}) for two reasons:\\
i) we assume that the density profile is continuous,\\
ii) it corresponds to a simple form of the energy equation.\\
The stationary associated quantities $\rho_*$ and $u_*$ are solution of:
$$\bs\ba{l}\rho_*(x)u_*(x)=-\rho_aV_a\cr
{d\over dx}(\rho_*(x)u_*(x)^2+p_*(x))=-\rho_*(x)g\cr
{d\over dx}(u_0(x)-L_0V_a{d\over dx}(Z(\rho_*(x))))=-V_aZ'(\rho_a\zeta_0)\rho_a\xi'(y_0)\delta_{x-y_0L_0}\ea\es$$
coming from
\be\label{KANLD}\bs\ba{l}\partial_t\rho+\mbox{div}(\rho {\vec u})=0\cr
\partial_t(\rho {\vec u})+\mbox{div}(\rho {\vec u}\otimes {\vec u} + pId)=\rho {\vec g}\cr
\mbox{div}({\vec u}+L_0V_a\nabla Z(\rho))=V_a\frac{1-\xi_0}{\xi_0}\delta_{x-y_0L_0}.\ea\es\ee
It is easy to see that the stationary solution is thus given by
$$u_0(y)=\bs\ba{l}-\frac{V_a}{\xi_0}, y\in ]-\infty, y_0]\cr
-\frac{V_a}{\xi(y)}, y\in ]y_0, +\infty[\ea\es$$
$$p_0(y)=\bs\ba{l}-\frac{\rho_aV_a^2}{\xi_0}-\rho_ag\xi_0L_0(y-y_0), y\in ]-\infty, y_0]\cr
-\frac{\rho_aV_a^2}{\xi(y)}-\rho_ag\int_{y_0}^y\xi(s)ds, y\in [y_0, +\infty[\ea\es$$
We have
\begin{prop}
The linear growth rate associated with the system (\ref{KANLD}) in the case $y_0$ fixed independant of $\alpha$ in the regime $\alpha\rightarrow 0$ and $\beta, \gamma$ fixed is

$$\gamma = \sqrt{\frac{1-\xi_0}{1+\xi_0}}-\frac{\beta}{\xi_0}$$
which corresponds to
$${\bar \gamma}=\sqrt{gk\frac{\rho_a - \rho_0}{\rho_a+\rho_0}} - kV_{blowoff}.$$
\label{quasivide}\end{prop}
We recall that for all $\xi_c$:

$$\ba{l}{E_0(\xi_c)\wedge E_{a, +}(\xi_c)\wedge F_+(\xi_c)\over \gamma\xi_c-\beta}={E_0\wedge E_{a, +}\over \gamma\xi_c - \beta}\wedge (i_4+\alpha{\tilde R}^0+{\alpha^2\over \lambda_+-\alpha\gamma\xi_c})\cr
= -\beta(\xi_cg_6^{\perp}+g_2^{\perp}+\xi_cg_5^{\perp}+g_1^{\perp})-(\gamma\xi_c+2\beta)g_3^{\perp}\cr
+\alpha[(\xi_c^2+\beta^2)f_1^{\perp}+\xi_c(\xi_c^2+(\beta+\gamma\xi_c)^2) f_4^{\perp}- \beta\gamma\xi_c^2(f_2^{\perp}+f_3^{\perp})]\cr
+{\alpha^2\over \lambda_+(\xi_c)-\alpha\gamma\xi_c}[\xi_c^2(\gamma\xi_c +\beta)(f_1^{\perp}+\xi_cf_4^{\perp})+ (\beta^2-\gamma^2\xi_c^2)(-\xi_c(\gamma\xi_c+\beta)f_4^{\perp}+\beta(f_1^{\perp}+\xi(f_2^{\perp}+f_3^{\perp})))]\ea$$
Note that the leading order term in $\alpha$ comes from the coefficient of $i_4$ and of ${\tilde T}^0$.
As the leading order term in $\alpha$ of the solution in the region $[y_0, +\infty[$ in the case where $y_0$ is independant of $\alpha$ is given by the leading order term of $f_-(1)\wedge E_{a, -}(1)$, we find $i_4\wedge E_{a, -}(1)=W_+$ as leading order term.
The Evans function is $$Ev(0, \beta,\gamma,\xi_0)=
\beta(\beta - \gamma\xi_0){\xi_0^{\nu}\over \beta^2-\gamma^2\xi_0^2}
[\beta (\xi_0-1)(T_2^0+T_3^0)+(2\beta + \gamma(\xi_0+1))T_1^0]$$
hence
$$Ev(0, \beta, \gamma,\xi_0)=\beta(\beta - \gamma\xi_0)^2{\xi_0^{\nu}\over \beta^2-\gamma^2\xi_0^2}
[\xi_0^2(1-\xi_0) - (1+\xi_0)(\beta + \gamma\xi_0)^2].$$
The theorem \ref{quasivide} is proven.
\paragraph{Remark}
This result is somewhat surprising: in fact when $\xi_0\rightarrow 0$ the limit of the growth rate disappears because it becomes negative: in the previous equality $\xi_0$ must verify
$$\xi_0^3+(\beta^2-1)\xi_0+\beta^2>0$$
to obtain a positive growth rate.
It is easy to check that this is never satisfied if $\beta^2>\mbox{max}_{\xi_0\in ]0, 1]} (\xi_0^2\frac{1-\xi_0}{1+\xi_0})$ and is satisfied only in a region $\xi_-(\beta), \xi_+(\beta)$ for $\beta^2<\mbox{max}_{\xi_0\in ]0, 1]} (\xi_0^2\frac{1-\xi_0}{1+\xi_0})$.

We may thus conclude that this model is not relevant according to the physical results.

 \section{Solution in the hypergeometric region}
\label{hypergeometrique} In this section, we identify the solution
when $y$ goes to $-\infty$. We show a result for $t\leq t_0$ finite.
As the coefficients behave as $|y|^{-{1\over \nu}}$, we call this
problem a Fuchsian problem, and, as the leading term of the equation
leads to a hypergeometrical equation, we call this the
hypergeometrical set-up. We introduce \be
\label{eta}\eta(t)={\xi(-{t\over \alpha\beta})\over \alpha^{1\over
\nu}}, p(t)=\alpha^{1\over \nu}{\gamma\over \beta}\eta(t).\ee Let
\be h(p)=\sum_{k=1}^p{\nu\over \nu -k}({\beta\over \gamma
}))^kp^{k-1}-\nu({\beta\over
\gamma})^{\nu}p^{\nu-1}[h_{\{\nu\}}({\beta p\over \gamma})+y_0].\ee
Note that $h$ is a complex valued function because $p$ is complex, but as $p\gamma^{-1}$ and $\gamma h(p)$ are real, $ph(p)$ is also real.
From lemma \ref{comportementxi}, we deduce \be {(\eta(t))^{\nu}\over
\beta} = {1\over \nu t}(1+ph(p)).\ee The system of equations on
$Z=\left(\ba{l}z_1\cr x_2\cr x_3\cr z_4\cr x_5\ea\right)$ obtained
from (\ref{systemfuchscompletinitial}) is \be\bs\ba{l}-{dz_1\over
dt}+{\eta^{\nu}\over \beta}(1-\alpha^{1\over \nu}\eta)
z_1+x_3+\alpha^{1\over \nu}{\gamma\over \beta^2}\eta^{\nu+1}z_4=0\cr
-{dx_2\over dt}+pz_1-x_3+\alpha^{2\over \nu}{1\over
\beta^2}\eta^{\nu+2}z_4=0\cr -{dx_3\over
dt}+2z_1+x_2-px_3+{\eta^{\nu}\over \beta}z_4=0\cr -{dz_4\over
dt}+z_1+{\eta^{\nu}\over \beta}z_4-x_5=0\cr -{dx_5\over
dt}+x_3-z_4=0.\ea\es\label{systemfuchscomplet}\ee Introduce the
matrices
$$M_0(p)=\left(\ba{ccccc}0&0&1&0&0\cr
p&0&-1&0&0\cr 2&1&-p&0&0\cr 1&0&0&0&-1\cr
0&0&1&-1&0\ea\right),NZ=\left(\ba{c}z_1\cr 0\cr z_4\cr z_4\cr
0\ea\right).$$ Using the equalities

$$\ba{l}{\eta^\nu\over \beta}={1\over \nu t} + {\alpha^{1\over \nu}\over (\nu t)^{1+{1\over \nu}}}{\gamma\over \beta}\beta^{1\over \nu}h(p)(1+ph(p))^{1\over \nu}\cr
\alpha^{1\over \nu}{\eta^{\nu+1}\over \beta}= {\alpha^{1\over \nu}\over (\nu t)^{1+{1\over \nu}}}(1+ph(p))^{1+{1\over \nu}}\beta^{{1\over \nu}-1}\cr
\alpha^{2\over \nu}{1\over \beta^2}\eta^{\nu+2}={p\over \gamma}\alpha^{1\over \nu}{\eta^{\nu+1}\over \beta}\ea$$
there exists a regular matrix ${\cal M}(p, \beta, \gamma)$ such that (\ref{systemfuchscomplet}) rewrites
$$-{dZ\over dt}+M_0(p)Z+{1\over \nu t}NZ+{\alpha^{1\over \nu}\over (\nu t )^{1+{1\over \nu}}}{\cal M}(p, \beta,\gamma)Z=0.$$
It is then natural to introduce

We associate with the system (\ref{systemfuchscomplet}) the model system

\be\label{systemfuchsmodel}{d{\tilde Z}\over dt}=M_0(0){\tilde Z}+{1\over \nu t}N{\tilde Z}\ee
and the model system of the vectorial product of three solutions of (\ref{systemfuchsmodel})
\be\label{systemmodel3}{d{\hat Z}\over dt}=M_0(0)^{(3)}{\hat Z}+{1\over \nu t}N^{(3)}{\hat Z}.\ee

We need the following

\begin{lemma}
\label{base0}
The matrix $M_0(p)$ has three eigenvalues $-p, 1, -1$. The eigenvector asociated with $\lambda_0(p)=-p$ is $e_0(p)=(1, -2, -p, 0, 1)$. The eigenvectors associated with the eigenvalue of multiplicity 2 $\lambda_-=-1$ are $e_-(p)=(1, -1-p, -1, 0, 1)$ and $f_-={1\over 2}(0,0,0,1,1)$. The eigenvectors associated with the eigenvalue of multiplicity 2 $\lambda_+=1$ are $e_+(p)=(1, p-1, 1, 0, 1)$ and $f_+={1\over 2}(0,0,0,1, -1)$.
\end{lemma}
We notice that $$Ne_0(p)=Ne_+(p)=Ne_-(p)=i_1=e_+(p)+e_-(p)-e_0(p)+f_+-f_-$$
and
$$Nf_+=Nf_-={1\over 2}(i_3+i_4)={1\over 4}(e_-(p)-e_+(p))+{1\over 2}(f_++f_-).$$

We also introduce the function

\be\psi_0(t)={\gamma\over \beta}\alpha^{1\over \nu}\int_{t_0}^t\eta(s)ds.\label{psi}\ee
which is the integral of $p$ and corresponds to the eigenvalue of largest real part of $-M_0^{(3)}$. This eigenvalue is associated with the eigenvector $P(p)$ such that $e_0(p)\wedge e_-(p)\wedge f_-=\frac{1-p}{2}P(p)$ and we find
\be\label{vecteurP}P(p)=(2+p)f_1^{\perp}+f_2^{\perp}+f_3^{\perp}-g_1^{\perp}-g_2^{\perp}-(2+p)g_3^{\perp}-g_5^{\perp}-g_6^{\perp}.\ee
The aim of this section is to prove the
\begin{theo}
\label{solutionfuchs}

1) There exists a unique solution $U_0(t, \alpha)$ of the system (\ref{systemfuchscomplet}) and $t_0>0$ such that there exists $C>0$ such that, for all $t\geq t_0$ we have

$$|U_0(t, \alpha)e^{2t+\psi_0(t, \alpha)}|\leq Ct^{1\over 2\nu}.$$

\be\label{limitet}\mbox{lim}_{t\rightarrow +\infty}U_0(t, \alpha)e^{2t+\psi_0(t, \alpha)}t^{-{1\over 2\nu}}=P(0),\ee
the vector $P(0)$ being given by (\ref{vecteurP}). Note that this condition is equivalent to (\ref{conditionunicite3}).

2) There exists a unique solution of the system (\ref{systemmodel3}) such that

$${\tilde Z}(t)e^{2t}t^{-{1\over 2\nu}}\rightarrow P(0)$$
when $t\rightarrow +\infty$.

3) For all $t_0>0$, there exists $C(t_0)>0$ such that for $t\geq
t_0$ the estimate \be\label{estimate**}|{\tilde
Z}(t)e^{2t}t^{-{1\over 2\nu}}- U_0(t, \alpha)e^{2t+\psi_0(t,
\alpha)}t^{-{1\over 2\nu}}|\leq C(t_0) \alpha^{1\over \nu}.\ee

\end{theo}

In a first paragraph, me make a reduction of the system to simplify its resolution.
\subsection{Reduction of the system}
The system (\ref{systemfuchscomplet}) rewrites

\be{dZ\over dt}=M_0(p)Z+{1\over \nu t}NZ + \alpha^{1\over \nu}t^{-1-{1\over \nu}}R(p)Z\ee
where $M_0(p)$ are analytic functions of $p$ for $p\leq p_0$ and $N$ is a constant matrix. The system satisfied by the vectorial product of three solutions of (\ref{systemfuchscomplet}) is

\be{dZ^{(3)}\over dt}=M_0(p)^{(3)}Z^{(3)}+{1\over \nu t}N^{(3)}Z^{(3)} + \alpha^{1\over \nu}t^{-1-{1\over \nu}}R(p)^{(3)}Z^{(3)}.\label{systemfuchs3complet}\ee
to which we associate its model system (\ref{systemmodel3}).

Introduce the eigenvalue of smallest real part of $M_0(p)^{(3)}$, (which is $\lambda_1(p)=-2-p$). Consider the matrix

$$\Lambda(p)=M_0(p)^{(3)}-\lambda_1(p)I$$
as well as the unknown $Y_*(t)=Y^{(3)}e^{-\int_{t_0}^t\lambda_1(p(s))ds}$. We obtain the system

\be{dY_*\over dt}=\Lambda(p)Y_*+{1\over \nu t}N^{(3)}Y_*+\alpha^{1\over \nu}t^{-1-{1\over \nu}}R(p)^{(3)}Y_*\quad (GS)\ee
associated with the model system
\be{d\over dt}({\tilde Y}e^{2t})=\Lambda(0)({\tilde Y}e^{2t})+{1\over \nu t}N^{(3)}({\tilde Y}e^{2t})\quad (MS)\ee
The aim of what follows is to identify the family of solutions of (GS) such that $Y_*(t)\simeq ct^ME_0$ when $t\rightarrow +\infty$, where $M={N_{11}\over \nu}$ and $E_0$ is the eigenvector associated with $\lambda_1(p)$.\\
\paragraph{Remark} The result that we obtain here depends heavily on the fact that the non zero eigenvalues of $\Lambda(p)$ denoted by $\lambda_i(p)$ satisfy
$$(\ln t)^{-1}\int_{t_0}^t\lambda_i(p)ds\rightarrow +\infty, t\rightarrow +\infty.$$

The first transformation uses $\Lambda(p)=(P(p))^{-1}D(p)P(p)$, where $P(p)$ is a transfer matrix and $D(p)$ is the matrix of eigenvalues of $\lambda(p)$, such that $D_{ii}(p)=\lambda_i(p)$, $\lambda_i(p)\leq \lambda_{i+1}(p)$ and all the eigenvalues are positive. We introduce $U(t, \alpha)=P(p)Y_*$. The system is

$$ {dU\over dt}=D(p)U+{1\over \nu t}P(p)[N^{(3)}+\nu ({\alpha\over t})^{1\over \nu}R(p)^{(3)}-P(p)^{-1}{dP\over dp}(\nu t{dp\over dt})]P(p)^{-1}U,$$
and using the relation

$${dp\over dt}=-\alpha^{1\over \nu}(\nu t)^{-1-{1\over \nu}}{\gamma\over \beta^{1-{1\over \nu}}}(1-{\beta\over \gamma}p)(1+ph(p))^{1+{1\over \nu}}$$
we end-up with the system

\be\label{SHG} {dU\over dt}=D(p)U+{1\over \nu t}M(t, \alpha)U.\ee
where the matrix $M$ writes

$$M(t, \alpha)=N^{(3)}+R(t, \alpha)\alpha^{1\over \nu}t^{-{1\over \nu}},$$
with the following estimate on $R$:

$$\exists \alpha^0, T_0>0, \forall \beta\in [\beta_0, {1\over \beta_0}], \forall \gamma, \Re\gamma\geq 0, |\gamma|\leq \beta_0^{-1}, |R(t, \alpha)|\leq C(\beta_0).$$
On each eigenspace of $D(p)$ we assume that $N$ is diagonal.

We denote by $E_i$ the eigenspaces of $D(p)$, $1\leq i \leq m$, $E_1$ and $E_m$ being of dimension 1. The unknowns $U$ are written $U_i\in K^{\mbox{dim}E_i}$.

The aim of the next paragraph is to construct iteratively the solution of the system (\ref{SHG}). We use the methods of Levinson \cite{Lev} and Hartmann \cite{Har}.
\subsection{Formal solution of the system}

It is necessary to begin with the computation of $U_m$, solution of

$$\ba{ll}{dU_m\over dt}=&\lambda_m(p)U_m+{1\over \nu t}(N_{mm}+\alpha^{1\over \nu}t^{-{1\over \nu}}R_{mm})U_m\cr
&+{1\over \nu t}\sum_{j=1}^{m-1}M_{mj}U_j.\ea$$
We consider the differential equation

$${dU_m\over dt}=\lambda_m(p)U_m+{1\over \nu t}(N_{mm}+\alpha^{1\over \nu}t^{-{1\over \nu}}R_{mm})U_m-f$$
 As we want to obtain a bounded solution of (\ref{SHG}), if we introduce $\phi_m(t)=\int_{t_0}^t\lambda_m(p(s))ds$,  this differential equation becomes

$${d\over dt}(U_m(t)e^{-\phi_m(t)}t^{-{N_{mm}\over \nu}})={1\over \nu t}\alpha^{1\over \nu}t^{-{1\over \nu}}R_{mm}(U_m(t)e^{-\phi_m(t)}t^{-{N_{mm}\over \nu}})-f(t)e^{-\phi_m(t)}t^{-{N_{mm}\over \nu}}.$$
If the function $U_m(t)e^{-\phi_m(t)}t^{-{N_{mm}\over \nu}}$ goes to a non zero finite limit when $t$ goes to $+\infty$, then $U_m(t)$ goes to infinity when $t$ goes to infinity {\bf under the sufficient condition $\frac{\Re \phi_m(t)}{\ln t}\rightarrow+\infty$ for $t\rightarrow +\infty$}, which is contradictory with the fact that we seek a bounded solution. Hence it is necessary (but not sufficient) that $U_m(t)e^{-\phi_m(t)}t^{-{N_{mm}\over \nu}}$ goes to 0 when $t$ goes to infinity.

The system rewrites

$${d\over dt}(U_m(t)e^{-\phi_m(t)}t^{-{N_{mm}\over \nu}}e^{\int_t^{+\infty}\alpha^{1\over \nu}s^{-1-{1\over \nu}}R_{mm}(s)ds})=- f(t)e^{-\phi_m(t)}t^{-{N_{mm}\over \nu}}e^{\int_t^{+\infty}\alpha^{1\over \nu}s^{-1-{1\over \nu}}R_{mm}(s)ds}.$$

We introduce the operator

\be T^{(m)}(f)(t)=e^{\phi_m(t)}t^{{N_{mm}\over \nu}}e^{-\int_t^{+\infty}\alpha^{1\over \nu}s^{-1-{1\over \nu}}R_{mm}(s)ds}\int_t^{+\infty}f(s)e^{-\phi_m(s)}s^{-{N_{mm}\over \nu}}e^{\int_s^{+\infty}\alpha^{1\over \nu}l^{-1-{1\over \nu}}R_{mm}(l)dl}.\ee
We verify
\be{d\over dt}(T^{(m)}(f))=\lambda_m(p)T^{(m)}(f))+{1\over \nu t}(N_{mm}+\alpha^{1\over \nu}t^{-{1\over \nu}}R_{mm})T^{(m)}(f))-f(t).\ee
Hence the equation on $U_m$ leads to the necessary relation

$$U_m(t)=-T^{(m)}[{1\over \nu t}\sum_{j=1}^{m-1}M_{mj}U_j]$$
which rewrites

\be U_m(t)={1\over \nu t}\sum_{j=1}^{m-1}T_j^{(m)}(U_j)\ee
with
$$T^{(m)}_j(U_j)=-\nu tT^{(m)}({1\over \nu t}M_{mj}U_j).$$
We replace this equality in the system satisfied by $(U_j)_{1\leq j\leq m}$. We obtain

\be 1\leq j\leq m-1 \quad{dU_j\over dt}=\lambda_j(p)+{1\over \nu t}\sum_{k=1}^{m-1}M_{jk}^{(1)}(U_k)\ee
where
\be M_{jk}^{(1)}(U_k)=M_{jk}(U_k)+{1\over \nu t}M_{jm}T^{(m)}_k(U_k).\ee
Of course, we notice that
\be M_{jk}^{(1)}(U_k)-n_{jk}U_k=(M_{jk}-N_{jk})U_k+{1\over \nu t}M_{jm}T^{(m)}_k(U_k)=O(\alpha^{1\over \nu}t^{-{1\over \nu}})\ee
hence its contribution is a regularizing operator.

The scheme of the proof is the same for all the terms of the vector $U$. At each stage, we obtain the system
\be 1\leq j\leq m-e\quad {dU_j\over dt}=\lambda_j(p)U_j+{1\over \nu t}\sum_{k=1}^{m-e}M_{jk}^{(e)}(U_k)\ee
where
\be M_{jk}^{(e)}(U_k)=M_{jk}^{(e-1)}(U_k)+{1\over \nu t}M_{jm-e+1}^{(e-1)}T^{(m-e+1)}_k(U_k).\ee
The operator $T^{(m-e)}$ of the following step is given by the solution $T^{(m-e)}(f)$ going to 0 at $+\infty$ of the equation

$${dU_{m_e}\over dt}=\lambda_{m-e}(p)U_{m-e}(t) + {1\over \nu t}M_{{m-e}{m-e}}^{(e-1)}(U_{m-e})(t)-f(t)$$
and the operators $T^{(m-e)}_k$, $1\leq k\leq m-e-1$ are given by

\be T^{(m-e)}_k(U_k)=-\nu t T^{(m-e)}[{1\over \nu t}M_{{m-e}k}^{(e)}(U_k)].\ee

The construction of the formal solution is done. We end up with the remaining equation on $U_1$, to which we cannot apply the previous method because the associated eigenvalue is 0. This system writes

\be {d\over dt}(t^{-{N_{11}\over \nu}}U_1(t))={1\over \nu t}t^{-{N_{11}\over \nu}}(M^{(m-1)}_{11}-N_{11})U_1.\ee
We deduce that there exists a constant $A$ such that

\be \label{equafinale}U_1(t)-At^{N_{11}\over \nu}=-t^{N_{11}\over \nu}\int_t^{+\infty}{1\over \nu s}s^{-{N_{11}\over \nu}}(M_{11}^{(m-1)}-N_{11})(U_1)(s)ds.\ee

\subsection{Proof of the convergence of the previous Volterra series}
The resolution ends up with the construction of the solution of (\ref{equafinale}). For this construction, it is necessary to study the regularity of all the operators $T^{(m-e)}$ (and of all the induced operators $T^{(m-e)}_k$ and $M^{(e)}_{jk}$).

For a given function $\psi$, such that $\psi(t)$ is increasing, going to $+\infty$ at $+\infty$, we introduce

$$\Lambda^K_{\psi}(t_0)=\{f\in C^{\infty}([t_0, +\infty[), \exists C, |f(t)|\leq Ct^Ke^{\psi(t)}\}.$$
We say that $\psi\in L_{m-e}^{\varepsilon_0}(t_0)$ if we have

$$\psi\in C^{\infty}([t_0, +\infty[), \psi'(t)-\Re\phi'_{m-e}(t)\leq -\varepsilon_0<0, t\geq t_0, t^2|\psi''(t)|\mbox{ bounded on }[t_0, +\infty[.$$

These notations being introduced, we consider the operator

$$K^{(1)}(f)=-t^{N_{11}\over \nu}\int_t^{+\infty}{1\over \nu s}s^{-{N_{11}\over \nu}}(M_{11}^{(m-1)}-N_{11})(f)(s)ds.$$
The equation that we intend to solve is

$$U_1-At^{N_{11}\over \nu}= t^{N_{11}\over \nu}K^{(1)}(U_1).$$
We notice that ${1\over \nu s}s^{-{N_{11}\over \nu}}(M_{11}^{(m-1)}-N_{11})(f)\in L^1([t_0, +\infty[)$ as soon as $f\in \Lambda_{N_{11}\over \nu}^0(t_0)$. Moreover, for $t_0^0$ given, there exists a constant $C_0$ such that

$$|{1\over \nu s}s^{-{N_{11}\over \nu}}(M_{11}^{(m-1)}-N_{11})(f)\leq C_0\alpha^{1\over \nu}s^{-1-{1\over \nu}}\mbox{max}_{l\in [s, +\infty[}(|f(l)l^{-{N_{11}\over \nu}}).$$
Hence we get the inequality

$$|A^{-1}K^{(1)}(At^{N_{11}\over \nu})|= |K^{(1)}(t^{N_{11}\over \nu})|\leq \alpha^{1\over \nu}C_0\int_t^{+\infty}s^{-1-{1\over \nu}}ds=C_0\nu \alpha^{1\over \nu}t^{-{1\over \nu}}.$$
Assume $\alpha\leq \alpha_0$ given. There exists a value of $t^0$, given by

$$t^0=\mbox{max}(t_0, (2C_0\nu)^{\nu}\alpha_0)$$
such that for $t\geq t^0$ we have

$$|K^{(1)}(t^{N_{11}\over \nu})|\leq {1\over 2}t^{N_{11}\over \nu}.$$
Hence, by induction, we get that

\be |(K^{(1)})^{(l)}(t^{N_{11}\over \nu})|\leq {1\over 2^l}t^{N_{11}\over \nu}, t\geq t^0\ee
from which we deduce the convergence of the series $\sum (K^{(1)})^{(l)}(t^{N_{11}\over \nu})$ and its bound by $2t^{N_{11}\over \nu}$. More precisely, we have, for all $\alpha\leq \alpha_0$
\be |(K^{(1)})^{(l)}(t^{N_{11}\over \nu})|\leq {1\over 2^l}({\alpha\over \alpha_0})^lt^{N_{11}\over \nu}, t\geq t^0\ee
hence the behavior when $\alpha\rightarrow 0$.

\section{The instability growth rate}
\label{finalO}
This section relies on the relation (\ref{evansfinal}) that we obtained in the third section. The scope of the present section is to derive a limit, when $\alpha\rightarrow 0$ of the Evans function $Ev(\alpha, \beta, \gamma)$. As the right hand side of (\ref{evansfinal}) depends on $t$, and the left hand side of (\ref{evansfinal}) is independant of $t$, we will study, for a given (suitable) $t>0$, the limit when $\alpha\rightarrow 0$ of the right hand side of (\ref{evansfinal}). This corresponds to the calculus of the limit when $\alpha\rightarrow 0$ of the functions $z_1, \xi(y)z_p$ $(p=2,3,4)$, $\xi(y) m_l$, $(l=1,2,3)$, $\xi^2(y)m_l$, $l=4,5,6$ and $R_j(t, \alpha)$, $L_k(t, \alpha)$ for $y=-\frac{t}{\alpha \beta}$. Recall that we restrict ourselves to the regime $\alpha\rightarrow 0$, $\beta$, $\gamma$ in a fixed compact set.\\
In all what follows, we introduce
\be\label{tauxreduit} r=\frac{\gamma}{\beta}.\ee
\subsection{Calculus in the overlapping region}
Introduce $t_*>0$ given such that $$\frac{\nu t_*}{\beta}<\frac{1}{2R}.$$
We introduce $\zeta(t, \alpha)=\zeta(\xi(-\frac{t}{\alpha\beta}), \alpha)=\frac{\alpha}{(\zeta(\xi(-\frac{t}{\alpha\beta}))^{\nu}}$. We have $$\xi(-\frac{t}{\alpha\beta})=\alpha^{\frac{1}{\nu}}(\zeta(t, \alpha))^{-\frac{1}{\nu}}.$$
The equation on $\xi$ gives $\zeta+O(\alpha^{\frac{1}{\nu}})=\frac{\nu t}{\beta}$, hence
\be\zeta(t, 0)=\frac{\nu t}{\beta}.\label{limitezeta}\ee
Hence we get (and it is the same for all other quantities)

$$\mbox{lim}_{\alpha\rightarrow 0}z_1(-\frac{t}{\alpha\beta}, \alpha)=\zeta(t, 0){\bar A}_1(\zeta(t, 0), 0).$$
We have, moreover

$$exp(-\int_0^{-\frac{t_0}{\alpha\beta}}\alpha\gamma\xi(y')dy')=exp(\alpha^{\frac{1}{\nu}}\int_0^{t_0}\frac{\gamma}{\beta}(\zeta(s, \alpha))^{\frac{1}{\nu}}ds)$$
hence the limit when $\alpha\rightarrow 0$ of this quantity is $1$.

For $\beta_0$ given and $\xi_0=\frac{1}{2}$, we identify $\alpha_0$ such that, for $\alpha\leq \alpha_0$ and $\zeta<\frac{1}{R}$, the asymptotic series defining $z_j, mÃ_p$ converges to an analytic function. As $\zeta(t_*, 0)=\frac{\nu t_*}{\beta}$, if one chooses $0<\frac{\nu t_*}{\beta}<\frac{1}{2R}$, that is $0<t_*<\frac{\beta}{2\nu R}$, there exists $\alpha_1>0$ such that, for $\alpha\leq \alpha_1$, $\zeta(t_*, \alpha)\leq \frac{3}{4R}<\frac{1}{R}$.\\
For $t\in [\frac{t_*}{2}, t_*]$, and $\alpha\leq \alpha_1$, $\zeta(t, \alpha)\leq \frac{3}{4R}$, hence the functions $(\xi(-\frac{t}{\alpha\beta}))^{\frac{\delta_p}{\nu}}z_p(-\frac{t}{\alpha\beta})$ and $(\xi(-\frac{t}{\alpha\beta}))^{\frac{d_q}{\nu}}m_q(-\frac{t}{\alpha\beta})$ are well defined through the expansion of section \ref{over}. Moreover, for $t\geq \frac{t_*}{2}$ there exists $C(\frac{t_*}{2})$ such that (\ref{estimate**}) holds. Hence for $t\in [\frac{t_*}{2}, t_*]$ we have the limit of the functions $R_p$ and $L_k$ when $\alpha\rightarrow 0$. We are now ready to prove our main Theorem:

\begin{theo}\label{lacanau}
Let $M$ be given. There exists $\alpha_*>0$ such that, for $0<\alpha<\alpha_*$, $\beta\in [\frac{1}{M}, M]$ and $|\gamma|\leq M, \Re\gamma\geq 0$, the Evans function of the system $Ev(\alpha, \beta, \gamma)$ does not vanish.
\end{theo}
\paragraph{Proof}
Recall that we proved that
$$Ev(\alpha, \beta, \gamma)=Ev_0(\beta, \gamma)+\alpha^{\frac{1}{\nu}}Ev_1(\beta, \gamma, \alpha).$$
The value of $Ev_0(\beta, \gamma)$ (which is expressed through (\ref{eva}), (\ref{zmfinal})) does not depend on $t_*$. Letting
the leading order term of $Ev_0(\beta, \gamma)$ go to 0 when
$t_*\rightarrow +\infty$ (which gives, of course, the value of $Ev_0(\beta, \gamma)$ because it does not depend on $t_*$) yields $r=1$ as only possible positive
solution for $Ev(0, \beta, \gamma)=0$ (see (\ref{egafina})). However, for $r=1$ the remaining leading order term in $t_*$ of $Ev(0,
\beta, \gamma)$ is not zero (see (\ref{egafina2})), hence a contradiction. \\
{\bf The system of Kull-Anisimov has no bounded complex growth rate $\gamma$
  when $Fr=O(\frac{1}{\varepsilon})$ in the limit $\varepsilon
  \rightarrow 0$.}
To be more precise, recall that
\be\label{eva}\ba{ll}Ev(\alpha, \beta, \gamma)=&e^{-(2+\frac{\mu(\alpha)}{\beta})t}t^{\frac{1}{2\nu}}\exp(-\int_0^{-\frac{t_0}{\alpha\beta}}\alpha\gamma\xi(y')dy')\times\cr
&[z_1(R_1-\frac{\xi}{\beta^2}L_1-L_2)+\xi z_2(R_2+\frac{\gamma}{\beta}L_1-L_4)\cr
&+\xi z_3(R_3+\frac{\gamma}{\beta}L_2-\frac{\xi}{\beta^2}L_4)\cr
&+\xi z_4(R_4-\frac{\gamma}{\beta}L_3-\frac{\xi}{\beta^2}L_5-L_6)\cr
&+\frac{\xi}{1-\xi}(m_1L_1+m_2L_2+m_3L_3)+\frac{\xi^2}{1-\xi}(m_4L_4+m_5L_5+m_6L_6)].\ea\ee
We introduce ${\bar A}_p(\zeta, \alpha), {\bar B}_q(\zeta, \alpha)$, $1\leq p\leq 4$, $1\leq q\leq 6$ the functions given by
\be \beta \zeta{\bar A}_p(\zeta, \alpha)=A_p(\zeta, \alpha), \beta\zeta{\bar B}_q(\zeta, \alpha)= B_q(\zeta, \alpha).\ee
We proved in Section \ref{over} that the solution $w_{+}^{(2)}$ of (\ref{s2}) which behaves as $e^{(\lambda_-(1)-\alpha\beta)y}$ when $y\rightarrow +\infty$, satisfying the condition (\ref{conditionunicite2}) is given through the relation (\ref{formulew}) $\sum_{j=1}^4z_jf_j+\sum_{p=1}^6 m_pg_p=e^{-\alpha\mu(\alpha)y}T^{(2)}w_+^{(2)}$,
 where $(z_j, m_p)$ are given by Proposition \ref{serieprolongee} through
\be\bs\ba{l}z_1= \beta + \beta\zeta{\bar A}_1(\zeta, \alpha)\cr
z_2=\gamma - \beta + \beta\frac{\zeta}{\xi}{\bar A}_2(\zeta, \alpha)\cr
 z_3= \beta + \beta\frac{\zeta}{\xi}{\bar A}_3(\zeta, \alpha)\cr
 z_4=- \beta +\beta \frac{\zeta}{\xi}{\bar A}_4(\zeta, \alpha)\cr
 m_1=\beta\frac{\zeta}{\xi}{\bar B}_1(\zeta, \alpha)\cr
 m_2= \beta\frac{\zeta}{\xi}{\bar B}_2(\zeta, \alpha)\cr
 m_3= \beta\frac{\zeta}{\xi}{\bar B}_3(\zeta, \alpha)\cr
 m_4= \beta\frac{\zeta}{\xi^2}{\bar B}_4(\zeta, \alpha)\cr
 m_5= \beta\frac{\zeta}{\xi^2}{\bar B}_5(\zeta, \alpha)\cr
 m_6= \beta\frac{\zeta}{\xi^2}{\bar B}_6(\zeta, \alpha).\ea\es\label{zmfinal}\ee
 Equality (\ref{evansfinal}) for a $t$ such that $\zeta(t, \alpha)<\frac{1}{R}$, along with $\xi=\alpha^{\frac{1}{\nu}}\zeta(t, \alpha)^{-{\frac{1}{\nu}}}$ yields

$$\ba{ll}Ev(\alpha, \beta, \gamma)&=e^{-(2+\frac{\mu(\alpha}{\beta})t}t^{\frac{1}{2\nu}}exp(-\int_0^{-\frac{t_0}{\alpha\beta}}\alpha\gamma\xi(y')dy')\frac{1-\xi(0)}{\xi(0)}\beta\cr
 &[(1 + \zeta{\bar A}_1)(R_1-\frac{\xi}{\beta^2}L_1-L_2)+(\xi (r-1)+\zeta{\bar A}_2)(R_2+\frac{\gamma}{\beta}L_1-L_4)\cr
&+(\xi  +\zeta{\bar A}_3)(R_3+\frac{\gamma}{\beta}L_2-\frac{\xi}{\beta^2}L_4)+(-\xi  +\zeta{\bar A}_4)(R_4-\frac{\gamma}{\beta}L_3-\frac{\xi}{\beta^2}L_5-L_6)\cr
&+\frac{\zeta}{1-\xi}[{\bar B}_1L_1+{\bar B}_2L_2+{\bar B}_3L_3+{\bar B}_4L_4+{\bar B}_5L_5+{\bar B}_6L_6].\ea$$

We proved in Section \ref{hypergeometrique} that the unique solution $U_0(t, \alpha)$ of the system (\ref{systemfuchscomplet}) satisfying the uniqueness condition (\ref{limitet})is
$$U_0(t, \alpha)=\sum_{p=1}^4f_p(t, \alpha)f_p^{\perp} + \sum_{q=1}^6 g_q(t, \alpha)g_q^{\perp}$$
and the functions $R_j$ and $L_p$ given by (\ref{RL}) satisfy the estimates
$$\sum_{j=1}^4|R_j(t, \alpha)-R_j^0(t)|+\sum_{p=1}^6|L_p(t, \alpha)-L_p^0(t)|\leq C(t_0)\alpha^{\frac{1}{\nu}}, t\geq t_0.$$
Note that the limit of the quantities $R_j^0$ (as well as $R_j$ when $t\rightarrow +\infty$) is known.
We now consider the equality on $Ev(\alpha, \beta, \gamma)$ when $\alpha\rightarrow 0$. The right hand side is independant of $t$ because the left hand side is independant of $t$. Hence its value can be considered at $t_*$. Once $t_*$ is fixed, we get the limit by taking $\alpha=0$, hence, using $\frac{\mu(0)}{\beta}=-r-1$

\be \label{131}\ba{ll}Ev(0, \beta, \gamma)&=e^{(r-1)t_*}t_*^{\frac{1}{2\nu}}\frac{1-\xi(0)}{\xi(0)}\beta\cr
 &[(1 + \zeta{\bar A}_1^0(\zeta))(R_1^0(t_*)-L_2^0(t_*))+\zeta{\bar A}_2^0(R_2^0+rL_1^0-L_4^0)\cr
&+\zeta{\bar A}_3^0(R_3^0+rL_2^0)+\zeta{\bar A}_4^0(R_4^0-rL_3^0-L_6^0)\cr
&+\zeta[{\bar B}_1^0L_1^0+{\bar B}_2^0L_2^0+{\bar B}_3^0L_3^0+{\bar B}_4^0L_4^0+{\bar B}_5^0L_5^0+{\bar B}_6^0L_6^0]]\ea\ee
where the relations are written at $t=t_*$ and at $\zeta = \zeta(t_*, 0)=\frac{\nu t_*}{\beta}$.\\
 Let $A_{p, j}$ and $B_{q, j}$ being given through
\be\beta \zeta {\bar A}_p^0(\zeta)=\beta\sum_{j=1}^{\infty}(\beta \zeta)^jA_{p,
j}, \beta \zeta {\bar B}_q^0(\zeta)=\beta\sum_{j=1}^{\infty}(\beta
\zeta)^jB_{q, j}.\ee By keeping only the leading order term in
$\xi(y)$ for each equation in the system (\ref{S1}) (which means
that we consider the order of each quantity $Z_p$ and $M_q$), we
obtain the recurrence system (\ref{recur}):

\be \bs\ba{l}\nu (j+1)A_{1, j+1}=A_{3, j}-rA_{4, j}-B_{3, j}\cr
(\nu (j+1)+1)A_{2, j+1}=rA_{2, j}-A_{3, j}-B_{1, j}-B_{5, j}\cr
(\nu (j+1)+1)A_{3, j+1}=A_{1, j}+A_{2, j}+rA_{3, j}-A_{4, j}-B_{2, j}-B_{6, j}\cr
(\nu (j+1)+1)A_{4, j+1}=-A_{3, j}+rA_{4, j}+B_{3, j}\cr
(\nu (j+1)+1)B_{1, j+1}=-rB_{1, j}-B_{2, j}-B_{4, j}-rB_{5, j}-A_{1, j}+(r^2-1)A_{2, j}\cr
(\nu (j+1)+1)B_{2, j+1}=B_{1, j}-rB_{2, j}+B_{3, j}-rB_{6, j}+rA_1{1, j}+(r^2-1)A_{3, j}\cr
(\nu (j+1)+1)B_{3, j+1}=B_{2, j}-rB_{3, j}+B_{6, j}+(1-r^2)A_{4, j}\cr
(\nu (j+1)+2)B_{4, j+1}=-B_{1, j}+B_{5, j}+rA_{2, j}+A_{3, j}\cr
(\nu (j+1)+2)B_{5, j+1}=B_{4, j}-B_{6, j}-A_{4, j}\cr
(\nu (j+1)+2)B_{6, j+1}=B_{3, j}+B_{5, j}+rA_{4, j}.\ea\es\label{recur}\ee
It is easy from this system to deduce that, under the hypothesis $|r|\leq M^{-2}$, there exists a constant $C>0$ such that $\sum_{p=1}^4|A_{p, j}|+\sum_{q=1}^6|B_{q, j}|\leq (\frac{C}{\nu})^j\frac{1}{j!}$, hence ensuring that the analytic expansions defining ${\bar A}^0_p$ and ${\bar B}^0_q$ are extendible for all $\zeta$. The study of this recurrence system is the aim of the next paragraph.
\subsection{Behavior of the equivalent solution}
In what follows, we study the recurrence system.\\
Consider $B^0(\xi, r)$ given by
$$\bs\ba{l}B_0(\xi, r)i_1=-ri_1+\frac{1}{\xi}i_3+i_4\cr
B_0(\xi, r)i_2=i_3\cr B_0(\xi, r)i_3=\xi i_1-i_2+i_5\cr B_0(\xi,
r)i_4=ri_4+(1-r^2)i_1-\frac{1}{\xi}i_2+\frac{r}{\xi}i_3\cr B_0(\xi,
r)i_5=i_3-\xi i_4+r\xi i_1.\ea\es$$ associated with the differential
equation $\frac{dY^0}{dt}=B_0(\xi, r)Y^0$. It is easy to check that
$B_0^{(2)}$ is associated with the system (\ref{S1reduit}), obtained
from (\ref{S1}) by taking into account the behavior of $Z_p$ and
$M_q$ that we obtained in Section \ref{over}:
\be\label{S1reduit}\bs\ba{l}{dZ_1\over dt}=\xi Z_3 -r\xi Z_4 -\xi
M_3\cr
  {dZ_2\over dt}=r Z_2 -Z_3  - (M_1+\xi M_5)\cr
  {dZ_3\over dt}= {1\over \xi}Z_1 +  Z_2+\gamma Z_3- Z_4 -(M_2+\xi M_6)\cr
  {dZ_4\over dt}=-  Z_3+r Z_4 +M_3\cr
  {dM_1\over dt}=-r M_1 -  M_2 - \xi ( M_4 +r M_5) -{1 \over \xi})Z_1 + (r^2-1)Z_2\cr
  {dM_2\over dt}=  M_1 -r M_2 +  M_3 -r\xi M_6 + {r\over \xi}Z_1 +(r^2-1)Z_3\cr
  {dM_3\over dt}= M_2 -r M_3 + \xi M_6 + (1-r^2)Z_4\cr
  {dM_4\over dt}=-{1\over \xi}M_1  +  M_5 +\frac{r}{\xi}Z_2 +\frac{1}{\xi} Z_3\cr
  {dM_5\over dt}=M_4 - M_6 -\frac{1}{\xi}Z_4\cr
  {dM_6\over dt}= {1\over \xi}M_3 +  M_5 +\frac{r}{\xi}Z_4.\ea\es\ee
We introduce
\be\label{source} F_0(y)=-\frac{r+1}{\nu}\xi^{-\nu}f_1
+\frac{1}{\nu+1}\xi^{-\nu-1}(f_3-g_1+rg_2).\ee
We have
\begin{lemma} Let $H$ be the unique solution going to zero when
$y\rightarrow 0$ of $$\frac{dH}{dy}=\varepsilon B_0^{(2)}(\xi(y),
r)H(y)+\varepsilon \frac{d}{dy}(F_0(y))$$ The coefficients $A_{p,
j}, B_{q, j}$ are the coefficients in the expansion in $\varepsilon
(\xi(y))^{-\nu}$ of $H$, where $\xi(y)^{-\nu}=\nu y$. The vector
$(\zeta {\bar A}^0_p, \zeta {\bar B}^0_q)$ is equal to
$H(\frac{\beta\zeta}{\nu \varepsilon})$ for $\zeta=\xi(y)^{-\nu}$.
\label{lemme9}
\end{lemma}

  We need to obtain the relation with the asymptotic expansion in
  $\alpha\beta$ of Section \ref{over}. For this purpose, we
  introduce $t=\varepsilon y$. We consider the following function, which
  will be tentatively a solution of (\ref{S1reduit})
\be\label{133}\bs\ba{l}Z_1(y)=\sum_{j=1}^{\infty}A'_{1,j}\varepsilon^j
(\xi(y))^{-\nu j}\cr Z_p(y)=\sum_{j\geq 1}A'_{p,
j}\varepsilon^j(\xi(y))^{-\nu j -1}, p=2, 3, 4\cr M_q(y)=\sum_{j\geq
1}B'_{q, j}\varepsilon^j(\xi(y))^{-\nu j -1}, q=1,2, 3\cr
M_l(y)=\sum_{j\geq 1}B'_{l, j}\varepsilon^j(\xi(y))^{-\nu j -2},
l=4,5,6.\ea\es\ee Consider for example the first equation of
(\ref{133}). It rewrites
$$\sum_{j=1}^{\infty}\varepsilon^j A'_{1, j}\frac{d}{dy}(\xi^{-\nu j})=
 \varepsilon (\sum_{j\geq 1}(A'_{3, j}-rA'_{4, j}-B'_{3, j})
 \xi^{-\nu j}\varepsilon^j)$$
hence
$$-\frac{d\xi}{dy}\sum_{j=0}\nu (j+1)A'_{1,
j+1}\varepsilon^{j+1}(\xi(y))^{-\nu (j+1)-1}=\sum_{j\geq
1}\varepsilon^{j+1}(\xi(y))^{-\nu j}(A'_{3, j}-rA'_{4, j}-B'_{3,
j}).$$ If we want to obtain $(A'_{p, j}, B'_{q, j})$ independant on
$\xi(y)$, it is a natural choice to write \be
-\frac{d\xi}{dy}=\xi^{\nu+1}\ee and, up to a change of origin in
$y$, we obtain \be\xi(y)^{-\nu}=\nu y.\label{reduxi}\ee Even with
this equation, there will still be an additional term in the
relation, related with $j=0$. The resulting equation on $Z_1(y)$
given in (\ref{133}) is
$$\frac{d}{dy}[Z_1(y)-\varepsilon A'_{1, 1}(\xi(y))^{-\nu}]=\xi Z_3 -r\xi Z_4 -\xi
M_3.$$ With these two relations, we obtain the recurrence relation
$$\nu(j+1)A'_{1, j+1}=(A'_{3, j}-rA'_{4, j}-B'_{3, j}), j\geq 1.$$
The same method applies to all the equations.\\
If we impose the initial conditions: \be\label{init}\ba{l}A'_{1,
1}=-\frac{r+1}{\nu}, A'_{2, 1}=A'_{4, 1}=0, A'_{3,
1}=\frac{1}{\nu+1},\cr
 B'_{1, 1}=-\frac{1}{\nu+1}, B'_{2,
1}=\frac{r}{\nu+1}, B'_{3, 1}=B'_{4, 1}=B'_{5, 1}=B'_{6, 1}=0\ea\ee
which corresponds to the source term $F_0(y)$ given by (\ref{source}), we
have the identity
$$A_{p, j}=A'_{p, j}, B_{q, j}=B'_{q, j}, \forall j\geq 1, \forall p=1, 2, 3, 4, \forall q=1, 2, 3, 4, 5, 6$$
where $(A'_{p, j},B'_{q, j})$ are the coefficients of the expansion
in $\varepsilon$ of the solution of (\ref{S1reduit}) with the source
term (\ref{source}) whereas $(A_{p, j}, B_{q, j})$ are the
coefficients of the expansion in $\zeta$ of $({\bar A}^0_p, {\bar
B}^0_q)$. Lemma \ref{lemme9} is proven.\\

\paragraph{Reduction of the problem}
Let us study the matrix $B_0$. Its eigenvalues are $0, 1$, and $-1$
and that associated eigenvectors are
  $$e_0=(\xi, -2, 0, 0, 1),
  e_1=(1-r, -\frac{1}{\xi}, 0, 1, 0), F_1=(\xi, -1, 1, 0, 1),$$
  $$e_{-1}=(r+1, -\frac{1}{\xi}, 0, -1, 0), F_{-1}=(\xi, -1, -1, 0, 1).$$
The inverse matrix is given by
  $$\bs\ba{l}2i_1=e_1+e_{-1}+\frac{1}{\xi}(F_1+F_{-1}-2e_0)\cr
  2i_2=F_1+F_{-1}-2e_0\cr
  2i_3=F_1-F_{-1}\cr
  2i_4=(1+r)e_1+(r-1)e_{-1}+\frac{r}{\xi}(F_1+F_{-1}-2e_0)\cr
  2i_5=F_1+F_{-1}-\xi (e_1+e_{-1}).\ea\es$$
 We rewrite a solution $V$ of $\frac{dV}{dy}=\varepsilon B_0V$
  as
  $$V=V_0e_0+V_1e_1+W_1F_1+V_{-1}e_{-1}+W_{-1}F_{-1}.$$
 We obtain
  $$\frac{dV}{dy}=\frac{dV_0}{dy}e_0+\frac{dV_1}{dy}e_1+\frac{dW_1}{dy}f_1
  +\frac{dV_{-1}}{dy}e_{-1}+\frac{dW_{-1}}{dy}f_{-1}+
  (V_0+W_1+W_{-1})\frac{d\xi}{dy}i_1+\frac{d\xi}{dy}\xi^{-2}(V_1+V_{-1})i_2,$$
  hence the associated system is
\be\bs\ba{l}\frac{dV_0}{dy}=\varepsilon\frac{d\xi}{dy}
[\xi^{-1}(V_0+W_1+W_{-1})+\xi^{-2}(V_1+V_{-1})]\cr
\frac{dV_1}{dy}=\varepsilon V_1 -\varepsilon
\frac12\frac{d\xi}{dy}(V_0+W_1+W_{-1})\cr
\frac{dW_1}{dy}=\varepsilon W_1-\varepsilon
\frac12\frac{d\xi}{dy}[\xi^{-1}(V_0+W_1+W_{-1})+
\xi^{-2}(V_1+V_{-1})]\cr \frac{dV_{-1}}{dy}=-\varepsilon
V_{-1}-\varepsilon \frac12\frac{d\xi}{dy}(V_0+W_1+W_{-1})\cr
\frac{dW_{-1}}{dy}=-\varepsilon W_{-1}-\varepsilon
\frac12\frac{d\xi}{dy}
[\xi^{-1}(V_0+W_1+W_{-1})+\xi^{-2}(V_1+V_{-1})].\label{Fsystem}\ea\es\ee
From the relation
  $$\ba{l}4(A_1f_1+A_2f_2+A_3f_3+A_4f_4+B_1g_1+B_2g_2
  +B_3g_3+B_4g_4+B_5g_5+B_6g_6)\cr
  =A_1(e_1+e_{-1}+\frac{1}{\xi}R_0)\wedge
  (e_1-e_{-1})+ (rA_2-B_1)R_0\wedge (e_1+e_{-1})+A_2R_0\wedge
  (e_1-e_{-1})\cr
  +(rA_3-B_2)(F_1-F_{-1})\wedge(e_1+e_{-1}+\frac{1}{\xi}R_0)
  +A_3(F_1-F_{-1})\wedge (e_1-e_{-1})+\cr
  (rA_4+B_3)(e_1+e_{-1}+\frac{1}{\xi}R_0)\wedge(F_1-F_{-1}-\xi(e_1+e_{-1}))\cr
+A_4(e_1+e_{-1})\wedge(F_1+F_{-1}-\xi(e_1+e_{-1}))\cr+B_4R_0\wedge(F_1-F_{-1})+
  B_5R_0\wedge
  (F_1+F_{-1}-\xi(e_1+e_{-1}))\cr +B_6(F_1-F_{-1})\wedge(F_1+F_{-1}-\xi(e_1+e_{-1})).\ea$$
 we deduce a new basis of $\Lambda^2(\RR^5)$ in which the
  coefficients are $A_1$, $rA_2-B_1$, $A_2$, $rA_3-B_2$, $A_3$,
  $rA_4+B_3$, $A_4$, $B_4$, $B_5$, $B_6$.
The source term (\ref{source}) in the basis of $\Lambda^2(\mathbb
R^5)$ associated with $e_0, e_{\pm 1}, F_{\pm 1}$ we find
\be\label{sourcebis}\ba{ll}F_0(y)=&-\frac{r+1}{\nu}\xi^{-\nu}\frac14
[(e_1+e_{-1})\wedge
(e_1-e_{-1})-\frac{1}{\xi}(e_1-e_{-1})\wedge (F_1+F_{-1}-2e_0)]\cr
&+\frac{1}{(\nu+1)}\xi^{-\nu-1}[(F_1-F_{-1})\wedge (e_1-e_{-1})
-\frac{1}{\xi}(e_1+e_{-1})\wedge (F_1+F_{-1}-2e_0)].\ea\ee 
\paragraph{End of the proof} The
coefficient of $e_1\wedge F_1$ in the source term is thus
$$-\frac14\xi^{-\nu-1}(\frac{r+1}{\nu}-\frac{1}{\nu+1}-\frac{1}{\xi}).$$
The theory of Fuchsian systems (see
Hartmann \cite{Har}) shows
  that there exists a constant $\alpha_*$ such that the projection of
  $H(y)$ on $e_1\wedge F_1$ behaves as $-\frac14t^{\alpha_*}e^{2t}(\frac{r+1}{\nu}-\frac{1}{\nu+1})e_1\wedge
  F_1$. The leading order term in $t_*$ of the Evans functions
$Ev(0, \beta, \gamma)$ writes
$$\ba{ll}Ev(0, \beta,
\gamma)&=e^{(r+1)t_*}t_*^{\frac{1}{2\nu}}t_*^{\alpha_*}\frac{1-\xi(0)}{\xi(0)}\beta\cr
 &[(R_1^0-L_2^0) + \zeta(t_*, 0)[{\tilde A}_1^0(R_1^0-L_2^0)+{\tilde A}_2^0(R_2^0+rL_1^0-L_4^0)\cr
&+{\tilde A}_3^0(R_3^0+rL_2^0)+{\tilde A}_4^0(R_4^0-rL_3^0-L_6^0)\cr
&+{\tilde B}_1^0L_1^0+{\tilde B}_2^0L_2^0+{\tilde
B}_3^0L_3^0+{\tilde B}_4^0L_4^0+{\tilde B}_5^0L_5^0+{\tilde
B}_6^0L_6^0]]\ea$$ where ${\tilde A}_p^0(t_*)= {\bar
A}^0_p(\zeta(t_*, 0))e^{-2t_*}t_*^{-\alpha_*}$, ${\tilde
B}_q^0(t_*)= {\bar B}^0_q(\zeta(t_*, 0))e^{-2t_*}t_*^{-\alpha_*}$.
Hence, as
$$-e_1\wedge F_1=-\xi
f_1+f_2-f_3+f_4+rg_1+(1-r)(g_2+g_3)-\frac{1}{\xi}(g_4+g_5)$$ and the
limit of $(R_p^0, L_q^0)$ when $t_*\rightarrow +\infty$ is given by
(\ref{conditionunicite3}), the limit of $$Ev(0, \beta,
\gamma)(e^{(r+1)t_*}t_*^{\frac{1}{2\nu}}t_*^{\alpha_*}\frac{1-\xi(0)}{\xi(0)}
\beta)^{-1}$$ is \be (r-1)(-\frac{r+1}{\nu}+\frac{1}{\nu+1}).\label{egafina}\ee As
we seek positive values of $r$, the factor $e^{(r+1)t_*}$ goes to
infinity when $t_*\rightarrow +\infty$, hence this limit is
necessarily 0. The two possible values of $r$ are thus $r=1$ and $r=-\frac{1}{\nu+1}$, hence
$r=1$.\\

We also check that, for $r=1$, the projection of $F_0(y)$ on the space associated
with the eigenvalue $+1$ of $B_0^{(2)}$, space generated by $e_1\wedge
e_0$ and $F_1\wedge e_0$, is $(\frac{1}{\nu}\xi^{-\nu -1}+\frac{1}{2(\nu+1)}\xi^{-\nu -2})e_1\wedge
e_0$. We have the relations

$$e_1\wedge e_0\wedge S=0, F_1\wedge e_0\wedge S=0$$ because
$e_1\wedge e_0=-\xi f_1+2f_2+f_4+g_1-\frac{1}{\xi}g_5$ and $F_1\wedge
e_0=-\xi g_1 -\xi g_2+2g_4+g_5+g_6$. The leading order term of the
projection on this space is 0.
We get $e_1\wedge
e_{-1}\wedge S=-4i_1\wedge i_2\wedge i_3\wedge i_4\wedge i_5$ because $e_1\wedge
e_{-1}=-2f_1+\frac{2}{\xi}f_2+\frac{2}{\xi}g_1$, and similarily $e_1\wedge
F_{-1}\wedge S=0$, $e_{-1}\wedge F_1
\wedge S=-2i_1\wedge i_2\wedge i_3\wedge i_4\wedge i_5$, $F_1\wedge
F_{-1}\wedge S=0$. As the vector $F_0(y)$ has the following projection
on the eigenspace associated with the eigenvalue 0
\be\label{egafina2}\ba{ll}Pr(F_0(y))&=-\frac{1}{2\nu}\xi^{-\nu}(2e_1\wedge
e_{-1}+\frac{1}{\xi}e_{-1}\wedge F_1-\frac{1}{\xi}e_1\wedge
F_{-1})\cr
&+\frac{1}{\nu+1}(e_{-1}\wedge F_1+\frac{1}{\xi}e_{-1}\wedge F_1
+ (1-\frac{1}{\xi})e_1\wedge F_{-1})\ea\ee
the associated leading order term gives a non-zero contribution, hence
a contradiction. The main theorem is proven.

\section{Annex: Volterra type expansions of solutions of differential systems}
\label{annexe1}
We consider a solution of (\ref{SL2}) associated with the growth rate $\lambda_-(1)-\alpha\beta$, which is the eigenvalue of smallest real part $-M_0^{(2)}(+\infty)$. We recall that there exists a regular function $\mu(\alpha)$ such that
\be \lambda_-(1)-\alpha\beta = -1+\alpha\mu(\alpha).\ee
We will prove that there exists a unique solution $w_+^{(2)}$ of (\ref{SL2}) satisfying (\ref{conditionunicite2}). In order to prove the existence and uniqueness of $w_+^{(2)}$ we prove the following

\begin{prop}
Let $U$ be the solution going to $(1, 0, ..., 0)$ at $+\infty$ of the model system:

$${dU_j\over dy}=\lambda_jU_j + (1-\xi(y))\sum_{k=1}^dN_{jk}(\xi(y))U_k(y), 1\leq j\leq d$$
where the properties of the complex numbers $\lambda_j$ and of the functions $N_{jk}$ are the following
$$|N_{jk}(\xi)|\leq M_0, \xi\geq\xi_0$$
$0=\lambda_1< \Re\lambda_2... \leq \Re\lambda_{d-1}$, $N_{dd}=0$.\\
The function $U$ is given by
$$U(y)=(1, 0, ...., 0)+(1-\xi(y))w(y).$$
\label{billet}
\end{prop}
Note that, in the hypothesis, $N_{dd}=0$ is only there for simplicity purposes and is obtained by considering in the last equation (in which we should have $(1-\xi)N_{dd}= \xi'\frac{N_{dd}}{\xi^{\nu+1}}$) the conjugation by the exponential of the primitive of  $\frac{N_{dd}}{\xi^{\nu+1}}$.

\paragraph{Construction by recurrence of the Volterra operators}
We prove Proposition \ref{billet} by recurrence. Consider the last equation (line $d$ of the previous system). We have
$$\frac{dU_d}{dy}=\lambda_dU_d+\sum_{k=1}^{d-1}(1-\xi(y))N_{dk}U_j(y).$$
This equation is equivalent to
$$\frac{d}{dy}(U_de^{-\lambda_d y})= \sum_{k=1}^{d-1}e^{-\lambda_d y}(1-\xi(y))N_{dk}U_j(y).$$
As $U_d$ is bounded, the limit of $U_de^{-\lambda_d y}$ is zero, otherwise $U_d$ would not be bounded. Integrating from $y$ to $y_0$ and letting $y_0$ go to infinity, we have
$$U_d(y)=\sum_{k=1}^{d-1}e^{\lambda_dy}\int_{y}^{+\infty}(\xi(s)-1)N_{dk}(\xi(s))e^{-\lambda_d s}U_j(s)ds,$$
and defining $K_{dj}^{(1)}$ through
\be\label{K1}K_{dj}^{(1)}(U)=(1-\xi(y))^{-1}e^{\lambda_dy}\int_{y}^{+\infty}(\xi(s)-1)N_{dj}(\xi(s))e^{-\lambda_d s}U(s)ds\ee
we obtain
\be\label{UD}U_d(y)=(1-\xi(y))\sum_{k=1}^{d-1}K_{dj}^{(1)}(U_j)(y).\ee
Let us study the properties of $K_{dj}^{(1)}$.\\
Consider $U$ such that, for $\xi\geq\xi_0$ (and $y_0$ such that $\xi_0=\xi(y_0)$) we have the estimate
\be\label{estimationN}\exists N, C, \forall y\geq y_0, |U(y)|\leq C(1-\xi(y))^N.\ee
We obtain easily the estimate (thanks to $\lambda_{d-p}\geq 0$)
$$\forall y\geq y_0, |K_{dj}^{(1)}(U)(y)|\leq \frac{M_0C}{(N+1)\xi_0^{\nu+1}}(1-\xi(y))^N.$$
Replacing (\ref{UD}) in the $d-1$ first equations of the differential system of Proposition \ref{billet}, we obtain
$$1\leq j\leq d-1, \frac{dU_j}{dy}=\lambda_jU_j+(1-\xi(y))\sum_{k=1}^{d-1}N_{jk}^{(1)}(U_k)$$
where
\be\label{N1}N_{jk}^{(1)}(U_k)=N_{jk}.U_k+N_{jd}K_{dj}^{(1)}.\ee
We get the estimate, under the assumption (\ref{estimationN}) on $U$
$$y\geq y_0, |N_{jk}^{(1)}(U)(y)|\leq M_0C(1+\frac{M_0}{(N+1)\xi_0^{\nu+1}})(1-\xi(y))^{N}.$$
This rewrites
$$y\geq y_0, |N_{jk}^{(1)}(U)(y)|\leq M_1C(1-\xi(y))^{N}$$
where $M_1=M_0(1+\frac{M_0}{(N+1)\xi_0^{\nu+1}})$.
The procedure proceeds as follows:\\
Let $p$ be an element of $0.. d-2$. We write the sequence of operators $K_{d-p, k}^{(p+1)}$, $1\leq k\leq d-p-1$, and $N_{jk}^{(p+1)}$, $j\leq d-p-1$, $k\leq d-p-1$ such that
\be\label{UP}U_{d-p}(y)=(1-\xi(y))\sum_{k=1}^{d-p-1}K_{d-p,k}^{(p+1)}(U_k)(y)\ee
\be\label{eqUJ}\frac{dU_j}{dy}=\lambda_jU_j+(1-\xi(y))\sum_{k=1}^{d-p-1}N_{j k}^{(p+1)}(U_k)(y)\ee
where \be\label{NP1}N_{j k}^{(p+1)}(U_k)=N_{j k}^{(p)}(U_k)+(1-\xi(y))^{-1}N_{j,d-p}^{(p)}((1-\xi)K_{d-p, k}^{(p+1)}(U_k))\ee
and the operators $K_{d-p, k}^{(p+1)}$ are constructed as follows using the auxiliary problem

\be\label{KP}\frac{dV}{dy}=\lambda_{d-p}V+(1-\xi(y))N_{d-p, d-p}^{(p)}(V)+f.\ee
Bounded solutions when $y\rightarrow +\infty$ of (\ref{KP}) are obtained, for $d-p\geq 2$, through
$$V(y)e^{-\lambda_{d-p}y}=\int_{y}^{+\infty}(\xi(s)-1)e^{-\lambda_{d-p}s}N_{d-p, d-p}^{(p)}(V)(s)ds-\int_y^{+\infty}f(s)e^{-\lambda_{d-p}s}ds.$$
This can be written, using $g(y)=V(y)e^{-\lambda_{d-p}y}$, under the form

\be\label{KP2}g=K_{d-p}(g)-\int_y^{+\infty}f(s)e^{-\lambda_{d-p}s}ds.\ee
Hence we may write, for $f$ satisfying (\ref{estimationN})
$$g= \sum_{l=0}^{\infty}(K_{d-p})^l(\int_{+\infty}^yf(s)e^{-\lambda_{d-p}s}ds)$$
which defines $T_{d-p}$ through
\be\label{TP}V=T_{d-p}(f).\ee
Replacing $f$ by $(1-\xi(y))\sum_{k=1}^{d-p}N_{d-p, k}^{(p)}(U_k)(y)$, we obtain the expression of $U_{d-p}$ in function of $U_k, 1\leq k\leq d-p-1$ as
$$U_{d-p}(y)=T_{d-p}((1-\xi)\sum_{k<d-p}N_{d-p,k}^{(p)}(U_k))$$
hence
\be\label{KDP}(1-\xi(y))K_{d-p, k}^{(p+1)}(U_k)= T_{d-p}((1-\xi)N_{d-p,k}^{(p)}(U_k)).\ee
\paragraph{Construction of the solution}
The last step of this recurrence is to construct the solution $U_1$. The equation on $U_1$ writes then
\be\label{U1}\frac{dU_1}{dy}=(1-\xi(y))N_{11}^{(d-1)}(U_1),\ee
where $N_{11}^{(d-1)}$ satisfy the relations of Proposition \ref{volterra00}. This is equivalent to
$$U_1(y)-U_1(y_0)=\int_{y_0}^y(1-\xi(s))N_{11}^{(d-1)}(U_1)(s)ds.$$
Using the limit $U_1(y_0)\rightarrow 1$ when $y_0\rightarrow +\infty$, we obtain
$$U_1(y)=1+\int_{+\infty}^y(1-\xi(s))N_{11}^{(d-1)}(U_1)(s)ds.$$
We introduce the operator $K^1$ such that $U_1=1+K_1(U_1)$. For $g$ satisfying the estimate (\ref{estimationN}), a consequence of proposition \ref{volterra00} is that
$$|K^1(g)(y)|\leq \frac{MC}{\xi_0^{\nu+1}(N+1)}(1-\xi(y))^{N+1}.$$
Introduce $U_1^0(y)=1$, and the sequence $U_1^N=K^1(U_1^{N-1})$. It is straightforward to show the inequality for $y\geq y_0$
$$|U_1^N(y)|\leq (\frac{M}{\xi_0^{\nu+1}})^N\frac{1}{N!}(1-\xi(y))^N$$
hence  as the series $\sum_{N=0}^{\infty}(\frac{M}{\xi_0^{\nu+1}})^N\frac{1}{N!}(1-\xi(y))^N$ is normally convergent for $y\geq y_0$, the series $\sum_{n=0}^{\infty}U_1^N(y)$ is normally convergent and is the only solution going to 1 as $y$ goes to $+\infty$ of (\ref{U1}). We have the estimate

$$|U_1(y)|\leq \exp(\frac{M}{\xi_0^{\nu+1}}).$$
which is of the form (\ref{estimationN}) for $N=0$.\\
\paragraph{End of the construction}
As $U_1$ satisfies an estimate of the form (\ref{estimationN}) for $N=0$, from the definition of the operators $K^{(p+1)}_{d-p, k}$ we deduce that
$$U_2(y)=(1-\xi(y))w_2(y)$$
where $|w_2(y)|\leq M \exp(\frac{M}{\xi_0^{\nu+1}})$. Replacing in the equality (\ref{UP}) for $p=d-3$ the inequality on $U_1$ and on $U_2$ we deduce that
$$U_3(y)\leq M (1-\xi(y))\exp(\frac{M}{\xi_0^{\nu+1}}).$$
By recurrence on $p$ we obtain the inequalities on $U_j$ for $j\not=1$, which proves Proposition \ref{billet}.
\paragraph{Estimates}
We prove in this Section the following
\begin{prop}
\label{volterra00}
The operators $N_{jk}^{(p)}$ and $K^{(p+1)}_{d-p, k}$ satisfy the following estimates for $g$ satisfying (\ref{estimationN}):
$$\forall y\geq y_0, |T(g)(y)|\leq \frac{CM}{\xi_0^{\nu+1}(N+1)}(1-\xi(y))^N.$$
\end{prop}
This is a consequence of the more precise proposition
\begin{prop}
\label{volterra000}
Under the assumption $g$ satisfy the estimate (\ref{estimationN}), the operators $N_{jk}^{(p)}$ satisfy the estimate

$$|N_{jk}^{(p)}(g)(y)|\leq CM_p(1-\xi(y))^N$$
where
$M_{p+1}=M_p(1+\frac{M_p}{\xi_0^{\nu+1}}\exp(\frac{M_p}{\xi_0^{\nu+1}}))$
and
$$|T_{d-p}(f)(y)|\leq \frac{C}{\xi_0^{\nu+1}}\exp(\frac{M_p}{\xi_0^{\nu+1}})(1-\xi(y))^N.$$
\end{prop}
We prove the second proposition by recurrence. The first estimate that we have to deduce from the recurrence assumption on $N_{jk}^{(p)}$ is the estimate on $K_{d-p, k}^{(p+1)}$.\\
To obtain this result, we have to study the behavior of $K_{d-p}$ through the estimate on $N_{jk}^{(p)}$. We have, for $f$ satisfying (\ref{estimationN}) for all $N$, the inequality
$$|K_{d-p}(f)(y)|\leq \int_y^{\infty}(1-\xi(s))^{N+1}Me^{\lambda_{d-p}(y-s)}ds.$$
As $N+1\geq 1$, we use $y-s\leq 0$ and $\lambda_{d-p}\geq 0$ to obtain
$$\forall y\geq y_0, |K_{d-p}(f)(y)|\leq M\int_y^{\infty}(1-\xi(s))^{N+1}ds\leq\frac{M}{\xi_0^{\nu+1}(N+1)}(1-\xi(y))^{N+1}.$$
Note that we increase the power of $(1-\xi(y))$ in the result.\\
The second estimate is based on the expression of $T_{d-p}$ obtained through (\ref{TP}). For $N\geq 1$ and $f$ satisfying (\ref{estimationN}), we have
$$|\int_y^{\infty}f(s)e^{-\lambda_{d-p}s}ds|\leq Ce^{-\lambda_{d-p}y}\frac{(1-\xi(y))^N}{N\xi_0^{\nu+1}}.$$
The estimate on $K_{d-p}$, deduced from the recurrence hypothesis is
$$|K_{d-p}(\int_y^{\infty}f(s)e^{-\lambda_{d-p}s}ds)|\leq C\int_y^{\infty}(1-\xi(s))^{N+1}M_pds\leq \frac{M_pC}{\xi_0^{\nu+1}}(1-\xi(y))^{N+1}.$$
Hence we obtain
$$|(K_{d-p})^l(\int_y^{\infty}f(s)e^{-\lambda_{d-p}s}ds)|\leq Ce^{-\lambda_{d-p}y}(\frac{M_p}{\xi_0^{\nu+1}})^l\frac{(1-\xi(y))^N}{N+l}\frac{1}{l!}$$
hence the series $\sum_{l\geq 0}e^{\lambda_{d-p}y}(K_{d-p})^l(\int_y^{\infty}f(s)e^{-\lambda_{d-p}s}ds)$ is normally convergent. It defines the function $T_{d-p}(f)$ and we have
$$|T_{d-p}(f)(y)|\leq C(1-\xi(y))^N\exp(\frac{M_p}{\xi_0^{\nu+1}}).$$
Finally, using this estimate on $T_{d-p}$ as well as the relation
(\ref{KDP}), for $f$ satisfying (\ref{estimationN}) we have $$|(1-\xi(y))K_{d-p, k}^{(p+1)}(f)|\leq C(1-\xi(y))^{N+1}\exp(\frac{1}{\xi_0^{\nu+1}})$$
 which gives
$$|K_{d-p, k}^{(p+1)}(f)|\leq C(1-\xi(y))^N\exp(\frac{M_p}{\xi_0^{\nu+1}}).$$
Using the estimate on $N_{j, d-p}^{(p)}$ we obtain
$$|(1-\xi)^{-1}N_{j, d-p}^{(p)}((1-\xi)K_{d-p, k}^{(p+1)}(f))|\leq M_pC(1-\xi(y))^N\exp(\frac{M_p}{\xi_0^{\nu+1}})$$
hence
$$|N_{jk}^{(p+1)}(f)(y)|\leq C(M_p+M_p\exp(\frac{M_p}{\xi_0^{\nu+1}})) (1-\xi(y))^N.$$
If we introduce $M_{p+1}=M_p(1+\frac{M_p}{\xi_0^{\nu+1}}\exp(\frac{M_p}{\xi_0^{\nu+1}}))$, we thus deduce
$$|N_{jk}^{(p+1)}(f)(y)|\leq CM_{p+1} (1-\xi(y))^N.$$
The proposition \ref{volterra000} is proven. As we used the assumption that $\Re\lambda_{d-p}>0$ to ensure that a solution of (\ref{KP}) is given through

$$\frac{d}{dy}(Ve^{-\lambda_{d-p}y})=e^{-\lambda_{d-p}y}(1-\xi(y))N_{d-p, d-p}^{(p)}(V)+fe^{-\lambda_{d-p}y}$$
hence
$$V(y)e^{-\lambda_{d-p}y}-V(y_*)e^{-\lambda_{d-p}y_*}= \int_{y_*}^ye^{-\lambda_{d-p}s}(1-\xi(s))N_{d-p, d-p}^{(p)}(V)(s)+f(s)e^{-\lambda_{d-p}s}ds.$$
If $Ve^{-\lambda_{d-p}y}$ has a limit, then this is in contradiction with the fact that $V$ is bounded for $\Re\lambda_{d-p}>0$, and this argument is no longer valid if $\lambda_{d-p}=0$. Hence this proves that the recurrence stops at $p$ such that $\lambda_{d-p}=0$ hence $p=d-1$.
The estimates of Proposition \ref{volterra00} are valid for $N_{jk}^{(p+1)}$, hence the recurrence proceeds till $\Re\lambda_{d-p}>0$. This recurrence processus stops for $\lambda_{d-p}=0$ because we cannot assert that the equation (\ref{KP}) has bounded solutions going to a constant for $y\rightarrow +\infty$.

\end{document}